\documentclass{amsart}

\usepackage[T1]{fontenc}
\usepackage[utf8]{inputenc}

\usepackage[english]{babel}

\usepackage[leqno]{amsmath}
\usepackage{amsthm}
\usepackage{amsfonts}
\usepackage{amssymb}

\usepackage{color}
\usepackage{wrapfig}

\usepackage{graphicx}

\usepackage{subfig}
\usepackage{caption}

\usepackage{pseudocode}

\usepackage{hyperref}
\hypersetup{
	unicode=true,
	colorlinks=true,
	urlcolor=blue,
	pdfauthor=Ewa Nowakowska,
	pdftitle=Overlap Measure for Gaussian Mixture Models,
	linkcolor=blue
}

\theoremstyle{plain}

\theoremstyle{remark}

\theoremstyle{definition}

\providecommand{\norm}[1]{\left\lVert#1\right\rVert}
\providecommand{\abs}[1]{\left\lvert#1\right\rvert}
\newcommand{\ud}{\mathrm{d}}
\newcommand{\pkg}[1]{{\fontseries{b}\selectfont #1}} 
 
\DeclareMathOperator{\diag}{diag}
\DeclareMathOperator{\rank}{rank}
\DeclareMathOperator{\argmin}{argmin}
\DeclareMathOperator{\argmax}{argmax}

\newenvironment{narrow}[2]{%
	\begin{list}{}{%
		\setlength{\topsep}{0pt}%
		\setlength{\leftmargin}{#1}%
		\setlength{\rightmargin}{#2}%
		\setlength{\listparindent}{\parindent}%
		\setlength{\itemindent}{\parindent}%
		\setlength{\parsep}{\parskip}}%
		\item[]}
	{\end{list}
}

\title[Overlap Measure for Gaussian Mixture Models]{Tractable Measure of Component Overlap \\ for Gaussian Mixture Models}

\thanks{Supported by National Science Center of Poland, DEC-2011/01/N/ST6/04174}

\begin{document}

\author{Ewa Nowakowska}
\author{Jacek Koronacki}
\author{Stan Lipovetsky}

\address{Ewa Nowakowska\\Institute of Computer Science\\Polish Academy of Sciences\\
		ul. Jana Kazimierza 5\\01-248 Warszawa\\Poland}
\email{ewa.nowakowska@ipipan.waw.pl}

\address{Jacek Koronacki\\Institute of Computer Science\\Polish Academy of Sciences\\
		ul. Jana Kazimierza 5\\01-248 Warszawa\\Poland}	
\email{jacek.koronacki@ipipan.waw.pl}

\address{Stan Lipovetsky\\GfK Custom Research North America\\Marketing \& Data Sciences\\8401 Golden Valley Rd.\\Minneapolis MN 55427\\USA}	
\email{stan.lipovetsky@gfk.com}



\begin{abstract}
	The ability to quantify distinctness of a cluster structure is fundamental for certain simulation studies, in particular for those comparing performance of different classification algorithms. The intrinsic integral measure based on the overlap of corresponding mixture components is often analytically intractable. This is also the case for Gaussian mixture models with unequal covariance matrices when space dimension $d > 1$. In this work we focus on Gaussian mixture models and at the sample level we assume the class assignments to be known. We derive a measure of component overlap based on eigenvalues of a generalized eigenproblem that represents Fisher's discriminant task. We explain rationale behind it and present simulation results that show how well it can reflect the behavior of the integral measure in its linear approximation. The analyzed coefficient possesses the advantage of being analytically tractable and numerically computable even in complex setups. 

\end{abstract}


\maketitle


\section{Introduction}\label{intro}


\subsection{Overview.} There are numerous measures designed to capture distance between distributions or -- more specifically -- overlap between components of a Gaussian mixture model. One of the oldest is the Bhattacharyya coefficient (see for instance \cite{bhatt} or \cite{fuku}), which reflects the amount of overlap between two statistical samples or distributions, a generalization of Mahalanobis distance described in \cite{day} or \cite{maclach}. In the context of information theory the most generic is the Kullback-Leibler divergence (see \cite{k-l}) -- a non-symmetric measure of difference between two distributions, also interpreted as expected discrimination information, which sets the link with possible classification performance. In \cite{inmana} an overlap coefficient is proposed that measures agreement between two distributions, it is applied to samples of data coming from normal distributions. Among more recent works, in \cite{dasgupta} a c-separation measure between multidimensional Gaussian distributions is defined, later developed in \cite{maitra} as exact-c-separation. In \cite{sun}, in the setup simplified to two clusters $k = 2$ and two dimensions $d = 2$, overlap rate is defined as a ratio of the joint density in its saddle point to its lower peak. The concept of ridge curve is introduced and further developed in \cite{ray} and \cite{ridge2}, generalized to arbitrary number of dimensions and clusters, turning the ridge curve into a ridgeline manifold of the dimension $k-1$.

All the measures use the parameters of the distributions to assess the overlap between the components and are typically formulated in terms of the underlying model. However, they can also be applied at the data level, as long as the class (or cluster) assignment is known. Then the model parameter estimates are used instead instead.


\subsection{Content.} In Section \ref{integral} we recall the generic concept of component overlap and its best linear approximation, we also show an example of an overlap assessment approach and point to common difficulties. Then, in Section \ref{fisher} we introduce what we refer to as Fisher's distinctness measure and we explain rationale behind it. Finally, in Section \ref{sim} we show results of a simulation study that illustrates how well the Fisher's coefficient can reflect the linear approximation of the original intractable overlap coefficient.


\section{Overlap of Distributions}\label{integral}


\subsection{Integral measure} 
\begin{wrapfigure}[13]{r}{190pt}
  	\begin{center}
  		\includegraphics[width=0.36\textwidth]{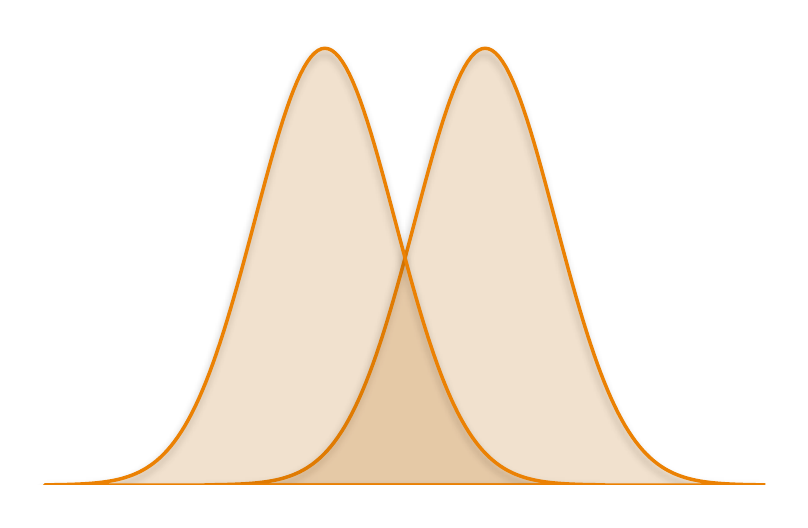}
  		\caption{Overlap (dark shadow) between $k = 2$ Gaussian components in $d = 1$ dimension, concept illustration.}
  		\label{overlap}
  	\end{center}
  \end{wrapfigure}

The most generic and natural coefficient of overlap between components is what follows directly from the mixture definition:
\begin{equation*}
	\text{MLE}_{\text{err}} = 1 - \int_{\mathbb{R}^d} \max \big( \pi_1 f_1 (\mu_1,\mathbf{\Sigma}_1), \ldots, \pi_k f_k (\mu_k,\mathbf{\Sigma}_k) \big ) (x) \ud x,
\end{equation*}
which for $k=2$ classes simplifies to
\begin{equation}\label{eq:47}
	\text{MLE}_{\text{err}} = \int_{\mathbb{R}^d} \min \big( \pi_1 f_1 (\mu_1,\mathbf{\Sigma}_1), \pi_2 f_2 (\mu_2,\mathbf{\Sigma}_2) \big ) (x) \ud x,
\end{equation} 
 where for $d \geq 1$ by $f_i$, $i = 1, \ldots, k$ we denote component densities and by $\pi_i$, $i = 1, \ldots, k$ their corresponding mixing factors. Throughout this work we will assume though that equal mixing factors are assigned to all the components, which corresponds to balanced cluster sizes at the sample level.  
 Coefficient \eqref{eq:47} measures the actual overlap between two probability distribution and for $d=1$ is illustrated in Figure \ref{overlap}. It coincides with intuitive understanding of components' overlap and with its expected behavior --- grows with increasing within cluster dispersion (or variance, for $d = 1$) and decreasing distance between cluster centers. Also, it exhibits a strong link with classification performance, setting the upper limit for possible predictive accuracy in terms of maximum likelihood estimation (MLE) (see for instance \cite{kmb}).
 Namely, best classification procedures based on likelihood ratio (MLE) or --- equivalently --- on its logarithm are given by 
 \begin{equation}\label{eq:45}
 	\text{loglik} (f_1, f_2) (x) = \text{log} \left(\frac{f_2(\mu_2, \mathbf{\Sigma}_2)(x)}{f_1(\mu_1, \mathbf{\Sigma}_1)(x)}\right).
 \end{equation}
 For the value of \eqref{eq:45} less than a constant observation $x$ is classified to the first cluster, to the second otherwise. Hence the area of overlap between the components, as given by \ref{eq:47}, corresponds to the expected proportion of observations that are incorrectly classified by MLE-classification rule, based on the (estimated) parameters of the mixture. Therefore \eqref{eq:47} is denoted by $\text{MLE}_{\text{err}}$ and alternatively referred to as $\text{MLE}$-misclassification or error rate. 

 The fundamental problem with formula \eqref{eq:47}, and also one of the reasons for numerous alternative approaches to overlap assessment, is that \eqref{eq:47} is hardly tractable for mixtures with different covariance matrices in higher dimensions. Handling it analytically would require integrating functions of Gaussian density over regions whose description often does not possess a tractable formulaic description either. Therefore, it can only be treated as a theoretical overlap coefficient for Gaussian mixture models and for practical applications replaced with other approaches.

 
 \subsection{Best linear approximation}\label{lin} The authors of \cite{ab} propose an approximation of \eqref{eq:47} --- best linear separator for $k=2$ Gaussian components in $d \geq 1$ dimensions and an algorithm to determine it for a given data set $X$. They derive a linear function of $x \in \mathbb{R}^d$ given by a vector $b \in \mathbb{R}^d$ such that for a given constant $c\in \mathbb{R}$ inequality $b^Tx \leq c$ classifies observation $x$ to the first cluster, while $b^Tx > c$ to the second. Vector $b$ and constant $c$ are obtained iteratively in order to minimize maximal probability of misclassification. As this approach will be used in our simulations, it is described below in more details following \cite{ab}. 
 
 For $x$ coming from component $l = 1,2$, $b^Tx$ has a univariate normal distribution with mean $b^T \mu_l$ and variance $b^T \mathbf{\Sigma}_l b$. As such, the probability of misclassifying observation $x$ when it comes from the first population $l=1$ equals
 \begin{equation}\label{eq:51a}
 	\mathbb{P}_1 \left(b^Tx > c \right) =  \mathbb{P}_1 \left(\frac{ b^Tx - b^T \mu_1}{b^T \mathbf{\Sigma}_1 b} > \frac{c -b^T \mu_1 }{b^T \mathbf{\Sigma}_1 b} \right) = 1 - \Phi \left( \frac{c -b^T \mu_1 }{b^T \mathbf{\Sigma}_1 b} \right) = 1 - \Phi \big( u_1 \big ),
 \end{equation}
 where $\Phi$ denotes cumulative distribution function for a univariate standardized normal distribution (centered at zero, with variance equal to one) and $u_1 = \frac{c -b^T \mu_1 }{b^T \mathbf{\Sigma}_1 b}$. Similarly, probability of misclassifying observation $x$ when it comes from the second population $l=2$ equals
 \begin{multline}\label{eq:51b}
 	\mathbb{P}_2 \left( b^Tx \leq c \right) =  \mathbb{P}_1 \left(\frac{b^Tx - b^T \mu_2}{b^T \mathbf{\Sigma}_2 b} \leq \frac{c -b^T \mu_2 }{b^T \mathbf{\Sigma}_2 b} \right) = \\
 	= \Phi \left(\frac{c -b^T \mu_2 }{b^T \mathbf{\Sigma}_2 b} \right) = 1 - \Phi \left(\frac{b^T \mu_2 -c}{b^T \mathbf{\Sigma}_2 b} \right) = 1 - \Phi \left( u_2\right),
 \end{multline}
 for $u_2 = \frac{b^T \mu_2 -c}{b^T \mathbf{\Sigma}_2 b}$. As $\Phi$ is monotonic, the task
 \begin{equation*}
      \max \big (\mathbb{P}_1(u_1), \mathbb{P}_2(u_2)\big)  \longrightarrow \min_{\substack{b \in \mathbb{R}^{d} \\ c \in \mathbb{R}}}
 \end{equation*}
 is equivalent to
 \begin{equation}\label{eq:48}
      \min (u_1, u_2)  \longrightarrow \max_{\substack{b \in \mathbb{R}^{d} \\ c \in \mathbb{R}}},
 \end{equation}
 which is more convenient to work with. As the objective is to find $b \in \mathbb{R}^d$ and $c \in \mathbb{R}$ that minimize maximal probability of misclassification, we will refer to the resulting procedure as a minimax procedure. Analytical formulation of admissible procedures for $b$ and $c$ leads to the following characterization
 \begin{equation}\label{eq:49b}
 	b = \left(t_1 \mathbf{\Sigma}_1 + t_2 \mathbf{\Sigma}_2\right)^{-1} (\mu_2 - \mu_1)
 \end{equation}
 and 
 \begin{equation}\label{eq:49c}
 	c = b^{T} \mu_1 + t_1 b^T \mathbf{\Sigma}_1 b = b^{T} \mu_2 - t_2 b^T \mathbf{\Sigma}_2 b, 
 \end{equation}
 where $t_1 \in \mathbb{R}$ and $t_2 \in \mathbb{R}$ are scalars. Minimax procedure is an admissible procedure with $u_1 = u_2$. As such, for $t = t_1 = (1-t_2)$ the following equality must hold
 \begin{equation}\label{eq:50}
 	0 = u_1^2 - u_2^2 = t^2 b^T \mathbf{\Sigma}_1 b - (1-t)^2 b^T \mathbf{\Sigma}_2 b = b^T \left[t^2 \mathbf{\Sigma}_1 - (1-t)^2 \mathbf{\Sigma}_2\right]b.
 \end{equation} 
 Equation \eqref{eq:50} for $t$ can be solved numerically by means of iterative procedure.
 
 With the above derivations, for a mixture of $k = 2$ components in $d \geq 1$ dimensions with parameters $\mu_1, \mathbf{\Sigma}_1$ and  $\mu_2, \mathbf{\Sigma}_2$ respectively, the following algorithm provides best linear separator in terms of minimizing the maximal probability of misclassification.
 
 \begin{center}
 	\begin{pseudocode}[ruled]{BestLinearSeparator}{\mu_1, \mathbf{\Sigma}_1, \mu_2, \mathbf{\Sigma}_2, prec}
 		$initialize $ incr,crit,t \\
 		\REPEAT 
 			$calculate $ b $ with \eqref{eq:49b}$ \\
 			$calculate $ crit $ with \eqref{eq:50}$ \\
 			\IF crit > prec \THEN t \GETS t - incr \\
 			\IF crit < -prec \THEN t \GETS t + incr \\
 			incr \GETS incr \cdot \frac{1}{2}
 		\UNTIL $criterion $ crit $ given by \eqref{eq:50} met with expected precision $ prec \\
 		$calculate $ c $ with \eqref{eq:49c}$ \\
 		$calculate $ u_1 $ and $ u_2 $ and the probabilities of misclassification with \eqref{eq:51a} and \eqref{eq:51b}$ \\
 		$calculate overall probability of misclassification $ \mathbb{P}_{\text{minmax}} = \max (\mathbb{P}_1 (u_1), \mathbb{P}_2 (u_2)) \\
 		\RETURN{\mathbb{P}_{\text{minmax}}, b, c, t}
 	\end{pseudocode}
 \end{center}
 
 Note that the value of assumed precision $prec$ must be given, while the values of scalar $t$, criterion $crit$ and increment $incr$ must be initialized. What is more, $\mathbb{P}_{\text{minmax}} = \mathbb{P}_1 (u_1) = \mathbb{P}_2 (u_2))$ as for the minimax procedure $u_1 = u_2$ must hold. Note, that $\mathbb{P}_{\text{minmax}}$ may be considered a measure of overlap as a result of linear approximation of criterion \eqref{eq:45}. If $\mathbf{\Sigma}_1 = \mathbf{\Sigma}_2$, formula \eqref{eq:45} and its linear approximation given by $b$ and $c$ coincide, which is sure not the case for $\mathbf{\Sigma}_1 \neq \mathbf{\Sigma}_2$.


 \subsection{The challenge of replacement.} Degree of overlap between mixture components is critical for classification performance and must be assessed for simulation purposes and comparison of classification methods, hence the interest in the topic. There are many measures proposed in the literature that possess the property of being tractable even in a complex setup, however it is highly required that their behavior reflects the behavior of $\text{MLE}_{\text{err}}$ based either on \eqref{eq:47} or on its linear approximation of the previous subsection. However, this is not always the case, as shown in the below example.

\textbf{E-distance.} The method for overlap assessment proposed in \cite{szekely} does not assume underling normal mixture model, however it can be very well applied in such setup. It is considered an extension to Ward's minimum variance method (see \cite{ward}) that formally takes both into account --- heterogeneity between groups and homogeneity within groups in data. For this purpose it uses joint between-within e-distance between clusters that constitutes the basis for agglomerative hierarchical clustering procedure the authors propose. They define e-\emph{distance} between two sets of observations $X_1 = \{x_{i_1} \colon c(i_1) = 1\}$, $n_1 = \abs{X_1}$ and $X_2 = \{x_{i_2} \colon c(i_2) = 2\}$, $n_2 = \abs{X_2}$ as 
\begin{multline}\label{eq:46}
	\mathrm{e}(X_1, X_2) = \frac{n_1 n_2}{n_1 + n_2} \left(\frac{2}{n_1 n_2} \sum_{i_1 \colon c(i_1) = 1} \sum_{i_2 \colon c(i_2) = 2} \norm{x_{i_1} -x_{i_2}} \right.+ \\
   \left.- \frac{1}{n_1^2} \sum_{i_1 \colon c(i_1) = 1}  \sum_{j_1 \colon c(j_1) = 1} \norm{x_{i_1} -x_{j_1}}    
   - \frac{1}{n_2^2} \sum_{i_2 \colon c(i_2) = 2}  \sum_{j_2 \colon c(j_2) = 2}\norm{x_{i_2} -x_{j_2}}\right).
\end{multline}
The value of e-distance between two resulting clusters may be considered a cluster structure distinctness measure. It is expected to reflect changes in within-cluster dispersion and between-cluster separation. It should also remain in tune with the theoretical structure distinctness measure given by likelihood ratio \eqref{eq:45} or its linear approximation from \cite{ab}. 
 
\begin{figure}[!h]
	\begin{narrow}{0.0\textwidth}{0.0\textwidth}
		\hfill
		\begin{minipage}[t]{0.49\textwidth}
			\centering
			\includegraphics[width=0.9\textwidth]{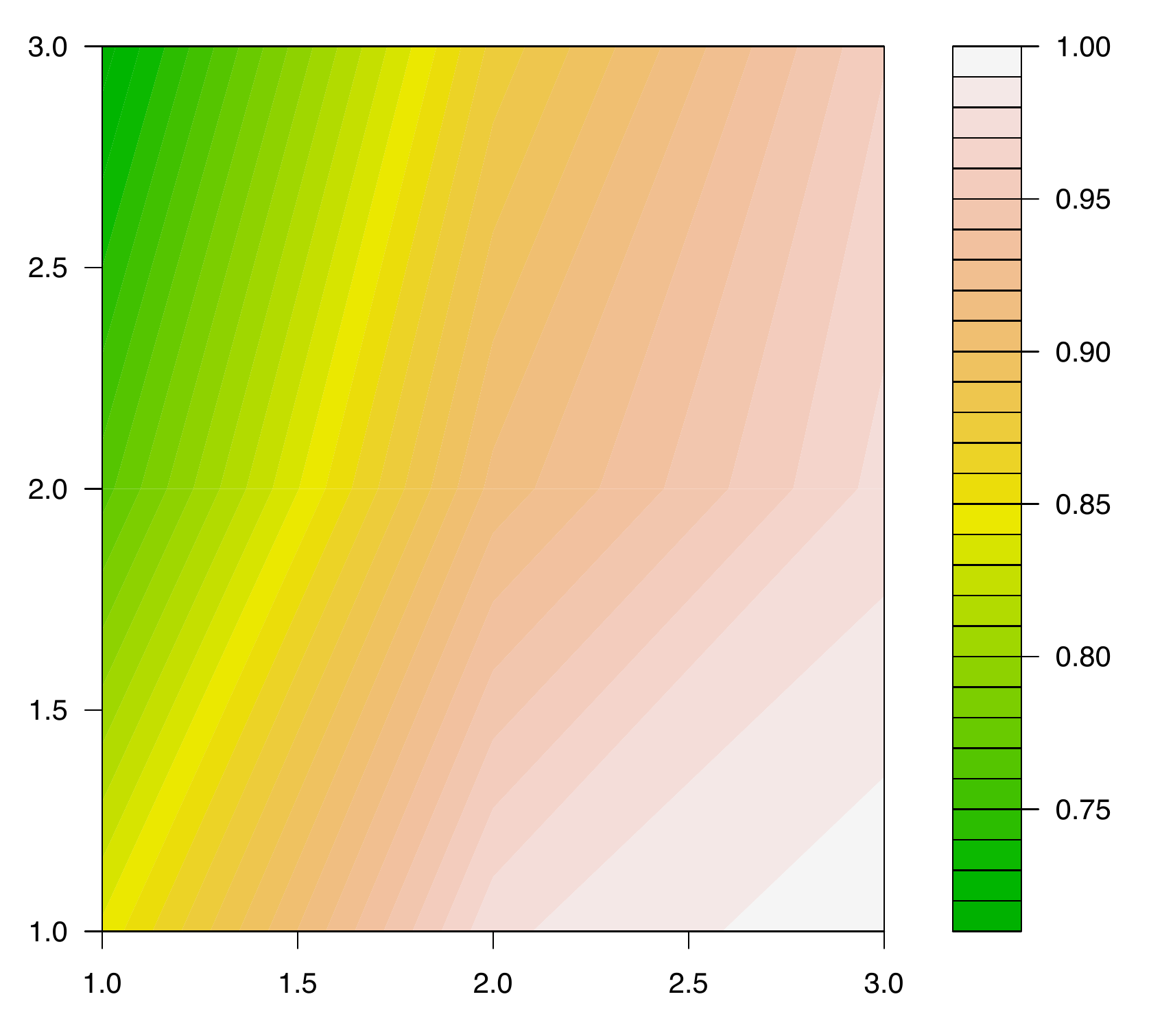}
			\captionof{figure}{Heatmap --- Anderson-Bahadur misclassification error w/r to growing between ($x$-axis) and within ($y$-axis) cluster dispersion.}
			\label{abclu}
		\end{minipage}
		\hfill
    	\begin{minipage}[t]{0.49\textwidth}
    		\centering
	    	\includegraphics[width=0.9\textwidth]{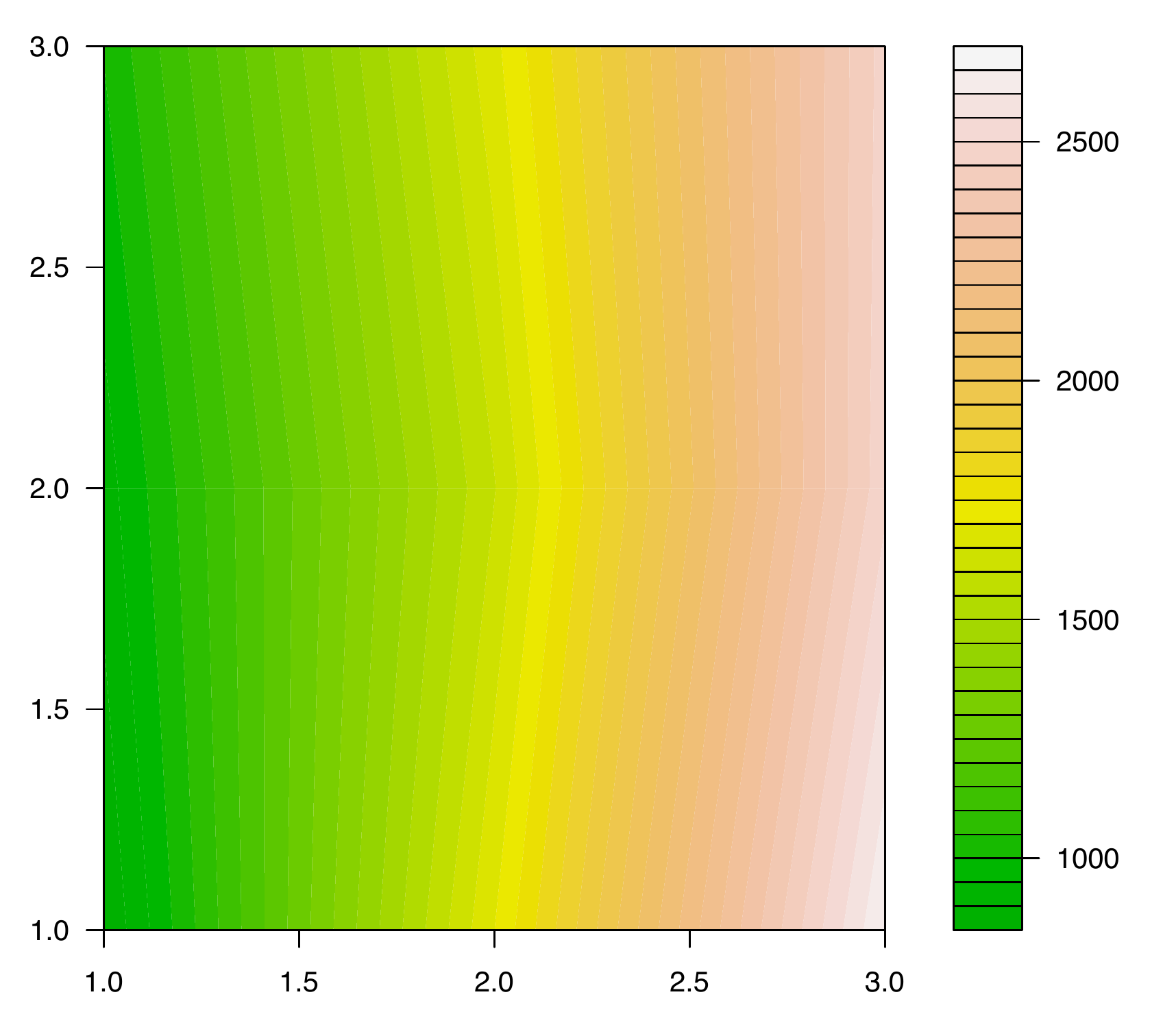}
    		\captionof{figure}{Heatmap --- Sz\'ekely-Rizzo e-distance \eqref{eq:46} w/r to growing between  ($x$-axis) and within ($y$-axis) cluster dispersion.}
    		\label{edist}
	    \end{minipage}
    	\hspace*{\fill}
    \end{narrow}
\end{figure}

Figures \ref{abclu} and \ref{edist} compare variability of structure distinctness measures based on Anderson-Bahadur (\cite{ab}) and Sz\'ekely-Rizzo (\cite{szekely}) proposals respectively. The former, similarly to the likelihood ratio theoretical distinctness measure,  does depend on both --- between-cluster distance and within-cluster dispersion, while the latter essentially remains insensitive to within cluster dispersion, depending entirely on the between class separation. This is an empirical insight which shows substantial discrepancy between behavior of theoretical and intuitive structure distinctness measure and e-distance given by \eqref{eq:46}, and hence points to another potential difficulty when trying to replace the integral coefficient.


\section{Fisher's distinctness measure}\label{fisher}


\subsection{Model and notation.} We consider a data set $X = ( x_1, \ldots, x_n)^T, \ X \in \mathbb{R}^{n \times d}$ of $n$ observations coming from a mixture of $k$ $d$-dimensional normal distributions 
\begin{equation*}
	f(x) = \pi_1 f_1(\mu_1, \mathbf{\Sigma}_1)(x) + \ldots + \pi_k f_k(\mu_k, \mathbf{\Sigma}_k)(x),
\end{equation*}
where
\begin{equation*}
	f_l(\mu_l,\mathbf{\Sigma}_l)(x) = \frac{1}{(\sqrt{2 \pi})^d \sqrt{\det\mathbf{\Sigma}_l}} e^{-\frac{1}{2}(x-\mu_l)^T
  \mathbf{\Sigma}_l^{-1}(x-\mu_l)}.
\end{equation*}
We call each $f_l(\mu_l, \mathbf{\Sigma}_l)$, $l = 1, \ldots, k$ a component of the mixture and each $\pi_l$, $l = 1, \ldots, k$ a mixing factor of the corresponding component (see \cite{kmb} or \cite{htf} and \cite{Lipovetsky:2013aa} or \cite{Lipovetsky:2012aa} for comparison with alternative approaches).  We assume that for all the components equal mixing factors are assigned $\pi_1 = \dots = \pi_k = \frac{1}{k}$. However, we allow different covariance matrices $\mathbf{\Sigma}_l$. Additionally, we assume large space dimension with respect to the number of components $d > k-1$ and take the number of components $k$ and class assignments as known. 

We use lower index to indicate data set when sample estimates of parameters are used. In particular, by $\mu_X \in \mathbb{R}^d$ we denote sample mean and by $\mathbf{\Sigma}_X \in \mathbb{R}^{d \times d}$ covariance matrix for a data set $X$. For notation ease we center the data at the origin $\mu_X = 0$. We assume the covariance matrix to be of full rank, $\rank(\mathbf{\Sigma}_X)  = d$. Let $T_X = n \mathbf{\Sigma}_X$ be the total scatter matrix for $X$. We recall that a simple calculation (see for instance \cite{kmb} or \cite{Lipovetsky:2013ab}) splits total scatter into its between and within cluster components $T_X = B_X + W_X$. By $\mu_{X,l}$ and $\mathbf{\Sigma}_{X,l}$ we denote empirical mean and covariance matrix for class $l$, where $l = 1, \ldots, k$. By $M_X = (\mu_{X,1}, \ldots, \mu_{X,k})$, $M_X \in \mathbb{R}^{d \times k}$ we understand a matrix of column vectors of cluster means. We assume the cluster means --- as a set of points --- to be linearly independent, so $\rank(M_X) = \min (d, k-1) = k-1$.


\subsection{Fisher's task as an eigenproblem.}\label{eigen} Originally (see \cite{fisher}), the separation was defined for $2$ classes in single dimension $v \in \mathbb{R}^d$ as the ratio of the variance between the classes to the variance within the classes
\begin{equation}\label{eq:1}
	F_o(v) = \frac{v^TB_Xv}{v^TW_Xv}.
\end{equation}
and then minimized over possible directions to find the linear subspace (Fisher's discriminant) that separates the classes best
\begin{equation*}
	v^* = \argmin (F_0(v)).
\end{equation*}
For our purposes we will use the formulation
\begin{equation}\label{eq8}
	F(v) = \frac{1}{1 + \frac{1}{F_o(v)}} = \frac{v^TB_Xv}{v^TT_Xv},
\end{equation}
which is equivalent to \eqref{eq:1} due to $T_X = B_X + W_X$ and yields the Fisher's subspace by maximizing over possible dimensions
\begin{equation}\label{eq:2}
	v^* = \argmax (F(v)).
\end{equation}

As multiplying $v$ by a constant does not change the result of  \eqref{eq:2}, it can alternatively be expressed as a constrained optimization problem, namely 
\begin{equation}\label{eq12}
	\begin{aligned}
		\max_{v \in \mathbb{R}^d } & & &v^T B_X v \\
		\text{subject to} & & &v^T T_X v = 1.
	\end{aligned}
\end{equation}
The corresponding Lagrange function defined as 
\begin{equation*}
	L(v; \lambda) =   v^T B_X v + \lambda\big(v^T T_X v - 1\big)
\end{equation*}
yields
\begin{equation*}
	\frac{\partial L(v; \lambda)}{\partial v} = 2 B_X v - 2 \lambda T_X v = 0,
\end{equation*}
so 
\begin{equation}\label{eq11}
	B_X v = \lambda T_X v
\end{equation}
must hold at the solution. Problem \eqref{eq11} is a generalized eigenproblem for two matrices $B_X$ and $T_X$. As we assume covariance matrix to be well-defined, total scatter matrix $T_X$ is invertible, however $T_X^{-1} B_X$ is not necessarily symmetric so it is a priori not obvious if the eigenvalues are real. Hence, a decomposition of the matrix $T_X$ is required to reduce the generalized eigenproblem to a standard eigenproblem for a transformed matrix. 

Solving a standard eigenproblem for $T_X$ we obtain 
\begin{equation}\label{eq13}
  T_X = A_{T_X} L_{T_X} A_{T_X}^{T}.
\end{equation} 
Note, that $A_{T_X}$ is orthonormal (i.e. $A_{T_X} A_{T_X}^{T} = \mathbf{I}$ so $ A_{T_X}^{-1} = A_{T_X}^{T}$). Replacing in \eqref{eq11} matrix $T_X$ with its spectral decomposition \eqref{eq13} we get
\begin{equation*}
	B_X v = \lambda A_{T_X} L_{T_X} A_{T_X}^{T} v = \lambda A_{T_X} L_{T_X}^{1/2} L_{T_X}^{1/2} A_{T_X}^{T} v,
\end{equation*}
then multiplying by $(A_{T_X} L_{T_X}^{1/2})^{-1}$ from the left and by $\mathbf{I}$ in the middle we transform it to
\begin{equation*}
	L_{T_X}^{-1/2} A_{T_X}^{T} B_X A_{T_X} L_{T_X}^{-1/2} L_{T_X}^{1/2} A_{T_X}^{T} v = \lambda L_{T_X}^{1/2} A_{T_X}^{T} v.
\end{equation*}
Now, substituting
\begin{equation*}
	\tilde{B} = L_{T_X}^{-1/2} A_{T_X}^{T} B_X A_{T_X} L_{T_X}^{-1/2} = \left(L_{T_X}^{-1/2} A_{T_X}^{T}\right) B_X \left(L_{T_X}^{-1/2} A_{T_X}^{T}\right)^{T}
\end{equation*}
and 
\begin{equation}\label{eq15}
	\tilde{v} = L_{T_X}^{1/2} A_{T_X}^{T} v
\end{equation}
we get a standard eigenproblem for $\tilde{B}$ 
\begin{equation}\label{eq14}
	\tilde{B} \tilde{v} = \lambda \tilde{v}.
\end{equation}
Solving \eqref{eq14} and using the inverse transformation of \eqref{eq15} 
\begin{equation}\label{eq16}
	v = A_{T_X} L_{T_X}^{-1/2} \tilde{v},
\end{equation}
we obtain the solution $v$ to the original problem \eqref{eq11}, corresponding to the same eigenvalue $\lambda$. In particular, it proves that with our model assumptions \eqref{eq11} can be reduced to a standard eigenproblem 
\begin{equation}\label{eq:36}
	T_X^{-1} B_X v = \lambda v,
\end{equation}
which takes the matrix form of
\begin{equation}\label{eq:43}
	\left ( T_X^{-1}B_X \right ) V = V L,
\end{equation}
where $L \in \mathbb{R}^{d \times d}$ is a diagonal matrix of eigenvalues in a non-decreasing order and $V \in \mathbb{R}^{d \times d}$ is a matrix of their corresponding column eigenvectors.

Note that there is another alternative formulation of the problem \eqref{eq11} via canonical correlation analysis (CCA), which may also come as a convenient way to see the task. In this setup Fisher's eigenvalues correspond to squared canonical correlation coefficients. We will not describe it here in details but we give references for interested readers. The approach, referred to as canonical discriminant analysis (CDA), was first mentioned in \cite{bartlett} and thoroughly described in \cite{dillon}. The overview of classical CCA is given for instance in \cite{kmb}.


\subsection{Motivation.} What we refer to as Fisher's distinctness measure was inspired by \cite{vemp}, where the idea of using the eigenproblem formulation of the Fisher's discrimination task and its respective eigenvalues for assessing certain properties of data was used. 

As explained in Subsection \ref{eigen}, Fisher's discriminant task can be stated in terms of eigenproblem given by \eqref{eq:43}. Then, its $(k-1)$ eigenvectors corresponding to the $(k-1)$ non-zero eigenvalues span the Fisher's subspace. Note that there are $k-1$ non-zero eigenvalues as according to the model assumptions $\rank(T_X) = d$ and $\rank(B_X)=k-1$ and $d > k-1$. Due to \eqref{eq:43} we have
\begin{equation*}
	V^T T_X^{-1} B_X V = L,
\end{equation*} 
so the eigenvalues capture variability in the spanning directions. As Fisher's task is scale invariant, the increase in variability may only be due to increase in between cluster scatter or decrease in within cluster scatter so it is expected to capture increase in structure distinctness very well. As squared canonical correlation coefficients (see references in Subsection \ref{eigen}), the eigenvalues remain in the interval of $[0,1]$ which also makes them easy to compare and interpret. Additionally, except for being easy to compute numerically, they are also convenient to handle analytically, so they can easily be used in simulations as well as formal derivations and justifications. What remains, is to propose function of the eigenvalues that could serve as structure distinctness coefficient and analyze its performance. This was done by means of simulation study and described in the next section.


\section{Simulation study}\label{sim}


\subsection{Overview.} Due to its analytical complexity \eqref{eq:47} is virtually intractable for mixtures with varied covariance matrices (heterogeneous) or of higher space dimension. However, it relatively easy undergoes simulations of Monte Carlo kind and can easily be approximated numerically with the best linear approximation described in subsection \ref{lin}. As such, it may be used as a reference measure and replaced with another coefficient that reflects its behavior but offers the advantage of being computable and analytically tractable, also in a more complex setup. 

The study was divided into two parts. In the first part two dimensional case was studied in details. Normal distribution was parametrized in a way that allowed for easy parameter control. Then all the possible combinations were tested and the influence of change in between cluster separation and within cluster dispersion was analysed. Three possible structure distinctness measures were compared --- exact integral measure \eqref{eq:47}, its best linear approximation described in subsection \ref{lin} and Fisher's eigenvalue. For two dimensional data, the maximum number of two clusters was analysed (due to the assumption of $d>k-1$), which led to one dimensional projections. Therefore, there was just single Fisher's eigenvalue to compare so the two dimensional step could not give grounds for function selection. The two dimensional study served as a thorough assessment of single Fisher's eigenvalue performance. 

In the second step multidimensional data was analyzed. Due to high number of possible mixture parameter combinations only a random selection was considered. This step was meant to confirm satisfactory performance of Fisher eigenvalues as input for structure distinctness measure. Higher dimensionality allowed for larger number of clusters, which resulted in $(k-1) > 1$ dimensionality of Fisher's subspace. As such, it also gave grounds for selecting appropriate function to transform $(k-1)$ eigenvalues into a single structure distinctness coefficient. Minimum $\lambda^X_{\text{min}}$ and average $\bar{\lambda}^X$ over Fisher's non-zero eigenvalues were calculated as follows 
\begin{equation}\label{eq:3a}
	\lambda^X_\text{min} = \min_{j \in \{1, \ldots, k-1 \} } \lambda_{j}^{T_X^{-1}B_X}
\end{equation}
and
\begin{equation}\label{eq:3}
	\bar{\lambda}^X = \frac{1}{k-1} \sum_{j = 1}^{k-1} \lambda_{j}^{T_X^{-1}B_X},
\end{equation}
and compared with the Monte Carlo estimates of the integral measure \eqref{eq:47}. Note that due to the larger number of classes allowed, wider comparisons with the best linear separator, defined for $k=2$ only, were infeasible. 

Note that although the original concept \eqref{eq:47} is defined in terms of overlap (similarity) between the components, what is naturally captured by either minimum or average over non-zero Fisher's eigenvalues, reflects the opposite behavior, so should rather be referred to as distinctness (dissimilarity) measure. Therefore we compare it with $(1 - \text{MLE}_{\text{err}})$ (or $(1-\mathbb{P}_{\text{minmax}})$), which is the probability of correct $\text{MLE}$ classification (or its best linear approximation). The transition from one to another is typically straightforward, however we point that out explicitly to avoid confusion or additional transformations of the coefficients.


\begin{center}
	\begin{pseudocode}[ruled]{TwoDimensionalDataGeneration}{r,\alpha,\lambda,q,k,N[]}
		\FOREACH $ cluster $ l \in \{1,\ldots,k\} \DO 
			\BEGIN
				\COMMENT{Determine cluster center $\mu$}\\
					\mu \GETS \left( r \cdot \sin \left((l-1) \cdot \frac{2 \pi}{k}\right), \ r \cdot \cos \left((l-1) \cdot \frac{2 \pi}{k}\right) \right) \\
				\COMMENT{Compute covariance matrix $\mathbf{\Sigma}$}\\
					D \GETS \diag(\lambda,q\cdot\lambda) \quad \COMMENT{dispersion and shape matrix} \\
					R \GETS \begin{pmatrix}
										\cos (\alpha) & -\sin (\alpha)  \\
										\sin (\alpha) & \cos (\alpha) 
									\end{pmatrix} \quad \COMMENT{rotation matrix} \\
					\mathbf{\Sigma} \GETS  R D R^T \\
				\COMMENT{Generate data} \\
					$draw $ N[l] $ observations$ \\
					$add cluster mean $ \mu $ to each observation$
			\END \\
		\RETURN{data}
	\end{pseudocode}
\end{center}

\subsection{Two-dimensional simulations.} To allow for easy control over mixture parameters, two dimensional mixture density was parametrized in a convenient way. Cluster centers were located on a circle around origin $(0,0)$ with radius $r$ that controlled between cluster distance. To allow for heterogeneity, for each cluster covariance matrix was determined separately. Within cluster dispersion was captured by the leading eigenvalue $\lambda = \lambda_1$, cluster shape by eigenvalues' ratio $q = \lambda_2 / \lambda_1$, and cluster rotation by rotation angle $\alpha$. Based on these parameters for each component mean vector and covariance matrix were computed. For each component the data was generated with the algorithm based on Cholesky decomposition, using affine transformation property for multivariate normal distribution. The detailed description of the algorithm is provided in \cite{gener}. Assuming the number of clusters is given by $k$ and $N \in \mathbf{R}^k$ contains desired cluster sizes, the above algorithm presents subsequent steps of data generation.

\begin{figure}[!ht]
	\begin{center}
		\includegraphics[scale=0.32]{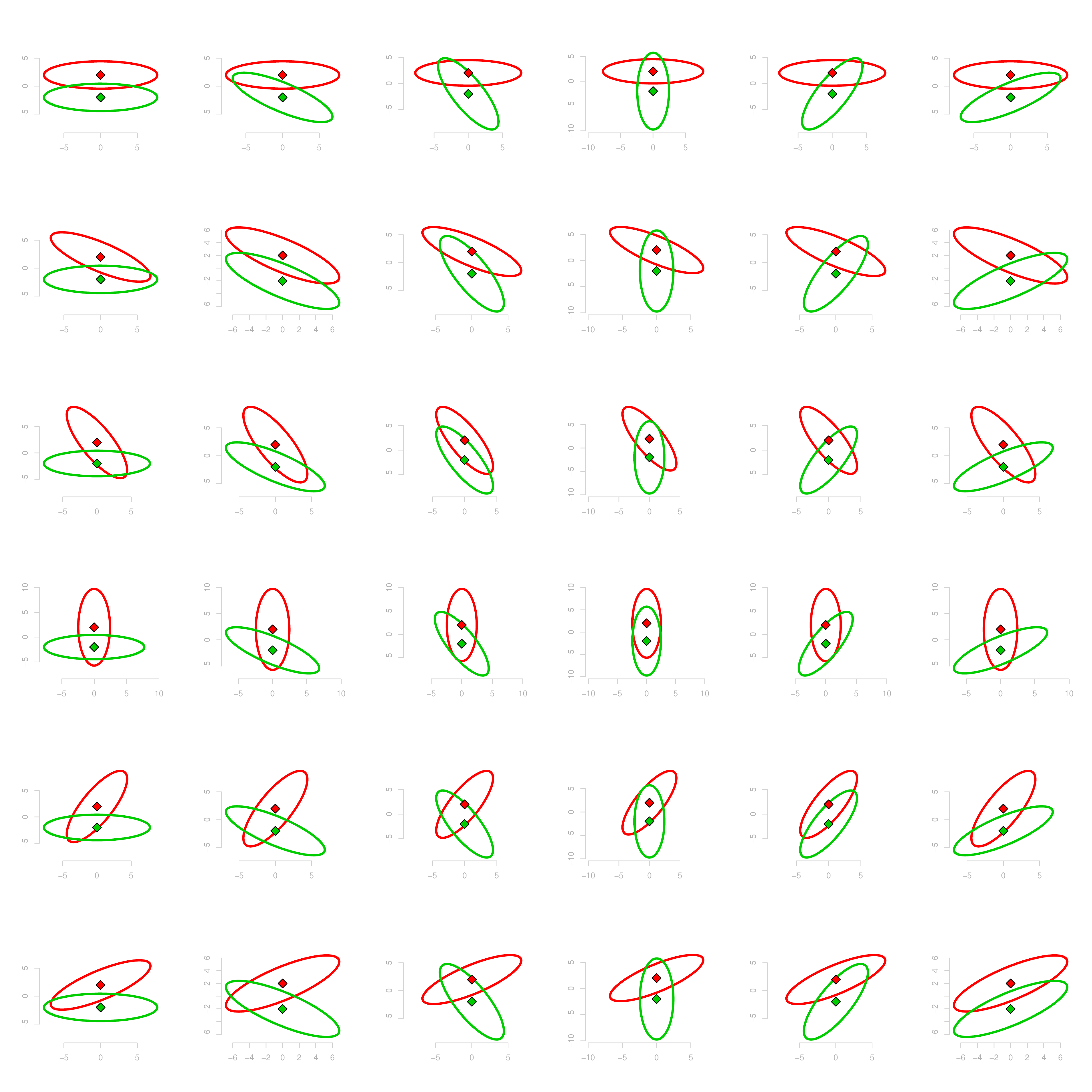}
		\captionof{figure}{Design of two dimensional simulations --- components' position with respect to each other.}\label{rys:11}
   \end{center}
\end{figure}

The simulation design is shown in Figure \ref{rys:11}, which presents all possible combinations of component position with respect to each other. Each of $i = 1, \ldots, 6$ rows corresponds to $i \cdot \pi/6$ angle rotation for the first (red) component, while each of $j = 1, \ldots, 6$ columns corresponds to $j \cdot \pi/6$ angle rotation for the second (green) component. Altogether it yields $36$ basic mixture positions. For each position an influence of a single factor is analyzed and this includes in particular -- increase in between cluster distance (Figures \ref{rys:12a} and \ref{rys:12b}), increase in within cluster dispersion for both (Figures \ref{rys:13a} and \ref{rys:13b}) and for first (Figures \ref{rys:14a} and \ref{rys:14b}) and second (Figures \ref{rys:15a} and \ref{rys:15b}) spanning direction only. The special case of spherical clusters is analyzed separately (Figures \ref{rys:16a} to \ref{rys:16d}). All the results are available in Appendix, Section \ref{app:2dim}.

\begin{figure}[!ht]
	\begin{center}
		\includegraphics[width=\linewidth]{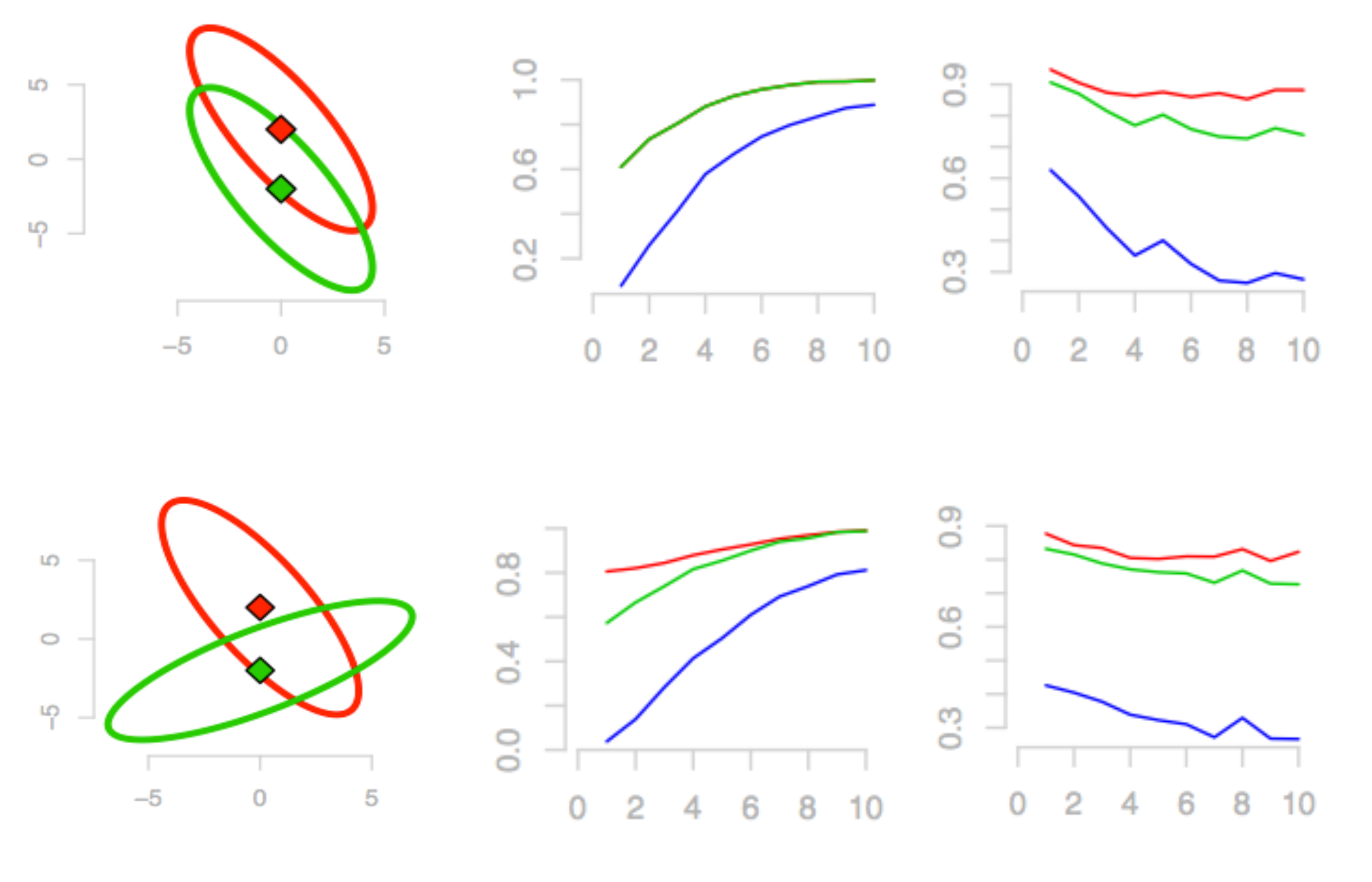}
        \caption{Impact of increasing between cluster distance (second column) and within cluster dispersion (third column) for mixtures in positions as indicated in the first column. Red line indicates exact (integral) structure distinctness, green --- its linear approximate and blue --- Fisher's eigenvalue.}
        \label{rys:20}
   \end{center}
\end{figure}

Example of what can be observed in all the charts is shown in Figure \ref{rys:20}. Even though the values for Fisher's eigenvalue are much smaller, their variability reflects behavior of the integral measure to a large extent. It is even more in tune with the linear estimate, which is to be expected given the linear nature of the Fisher's discrimination task. Note, that the best linear approximate gives the upper bound on the precision with which a linear concept may reflect behavior of the non-linear integral measure. Also, it gives upper limit on classification accuracy using linear classifiers, which is the case of Fisher discriminant. Note also, that the component position in the upper row indicates homogeneity (i.e. equal covariance matrices for both components). This property is lost when within cluster variability increases for one of the components (last column). However, it remains when only between cluster distance is affected (middle column). Therefore, exact integral measure and its linear estimate overlap in this case.


\subsection{Multi-dimensional simulations.}
 In higher dimensions direct analytical control over distance and dispersion of mixture parameters is much more complex. Additionally, there are many more combinations to examine. As such, the simulations were reduced to randomly chosen mixture parameters' combinations corresponding to the mixture position. For each position the impact of increasing between cluster distance and within cluster dispersion was analysed. The study was designed to verify adequacy of the information carried by the Fisher's eigenvalues and to select its appropriate function to serve as the structure distinctness coefficient. Results are attached in Appendix \ref{app:2dim} in Figures \ref{rys:17a} to \ref{rys:17d}. In each row charts for random but fixed set of cluster means are presented. Similarly, the set of covariance matrices is random but fixed in each column. Mean vectors and covariance matrices in $d$ dimensions were determined using \textbf{R} package \pkg{clusterGeneration}, which implements the ideas described in \cite{joe} and \cite{kurowicka}. Additionally, mean coordinates are re-scaled to lie in the interval $[-3\sqrt{d},3 \sqrt{d}]$ which corresponds to the range of the maximum three standard deviations for covariance matrix. As such, the possible overlap between components stretches from complete to negligible. 
 
 \begin{figure}[!ht]
 	\begin{center}      
        	 \includegraphics[width=\linewidth]{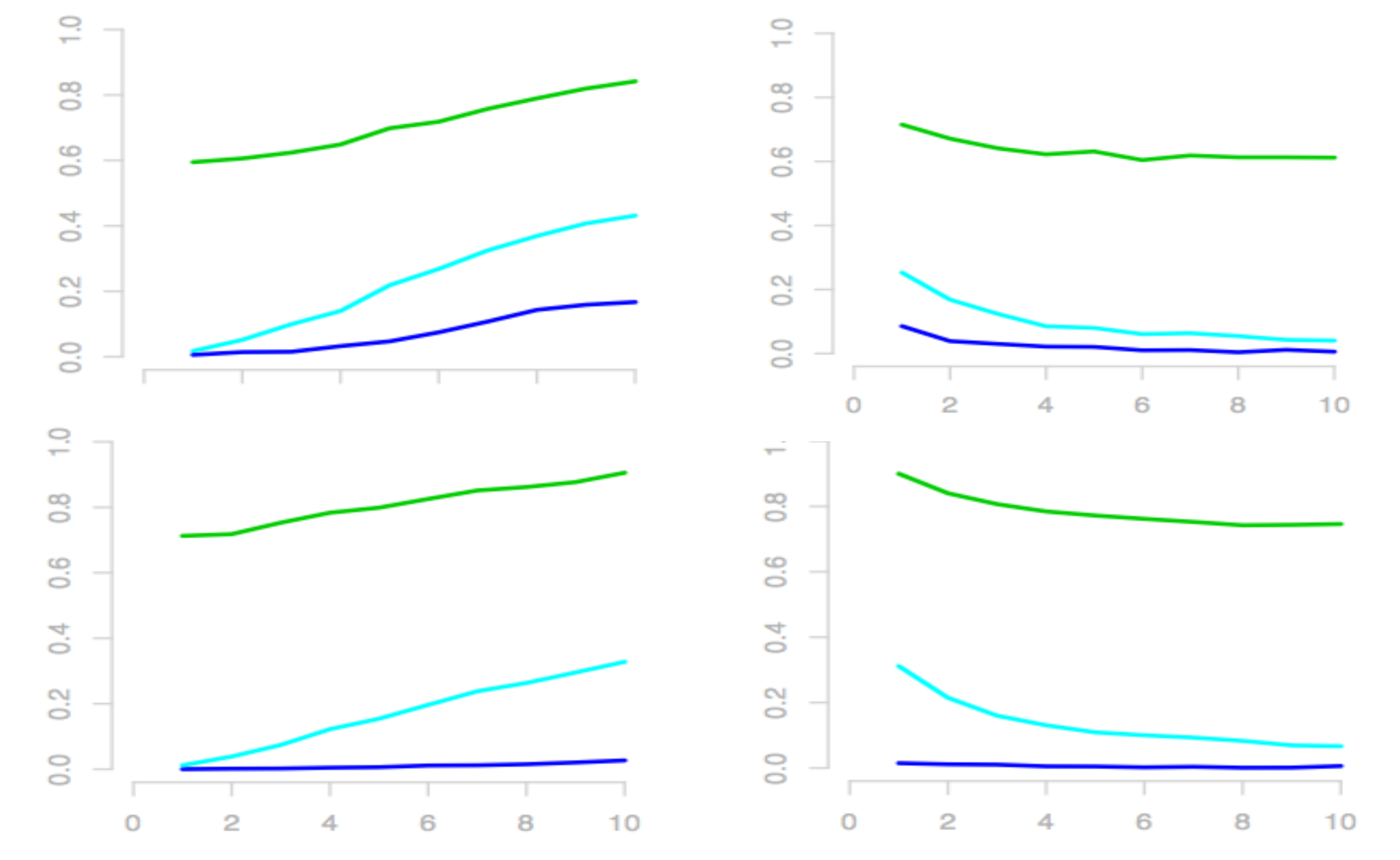}
         \caption{ Effect of increasing between cluster distance (left column) and within cluster dispersion (right column). Upper row gives results for three dimensional simulations, while bottom row for five dimensional case. Green line indicates Monte Carlo estimate of the integral structure distinctness, turquoise average non-zero Fisher's eigenvalue, while blue --- Fisher's smallest non-zero eigenvalue.}
         \label{rys:21}
    \end{center}
 \end{figure}
 
 Again, what can be observed in all the simulation plots in Appendix \ref{app:2dim} is illustrated in Figure \ref{rys:21}. Behavior of average Fisher's eigenvalue as given by \eqref{eq:3} reflects variability of the integral measure. At the same time, minimum non-zero Fisher's eigenvalue \eqref{eq:3} is less sensitive and therefore captures the changes in distinctness to a lesser extent, which becomes even more apparent as the number of dimensions increases. As such, the average non-zero Fisher's eigenvalue tends to outperform the minimum non-zero Fisher's eigenvalue and therefore the former shall be recommended as the distinctness coefficient.


\section{Conclusions.}

In this work we derive and motivate measure of distinctness (or alternatively -- overlap) between clusters of data, generated from a Gaussian mixture model. The approach uses alternative formulation of Fisher's discrimination task, which is stated in terms of a generalized eigenproblem. We show the task is well posed in the context of the assumed model and can be reduced to a standard eigenproblem with real eigenvalues. We then express the distinctness coefficient as the average eigenvalue over the non-zero eigenvalues of the solution. We compare the behavior of the coefficient with the generic (integral) measure of structure distinctness defined in terms of the actual overlap between the corresponding distributions and its best linear approximation. Although the values of the Fisher's coefficient are lower than the values of actual overlap, their dynamic reflects very well the behavior of the generic integral measure and even better -- its best linear approximation. As opposed to the generic integral measure and its best linear approximation, the Fisher's coefficient offers the advantage of being not only numerically easily computable but also analytically tractable, even in a complex setup, regardless of the dimensionality of the space and heterogeneity of covariance matrices.


\bibliographystyle{elsarticle-num}
\bibliography{mybibliography}


\newpage
\appendix


\section{Simulation results}\label{app:2dim}

\begin{figure}[!ht]
	\begin{center}      
       	 \includegraphics[width=\textwidth]{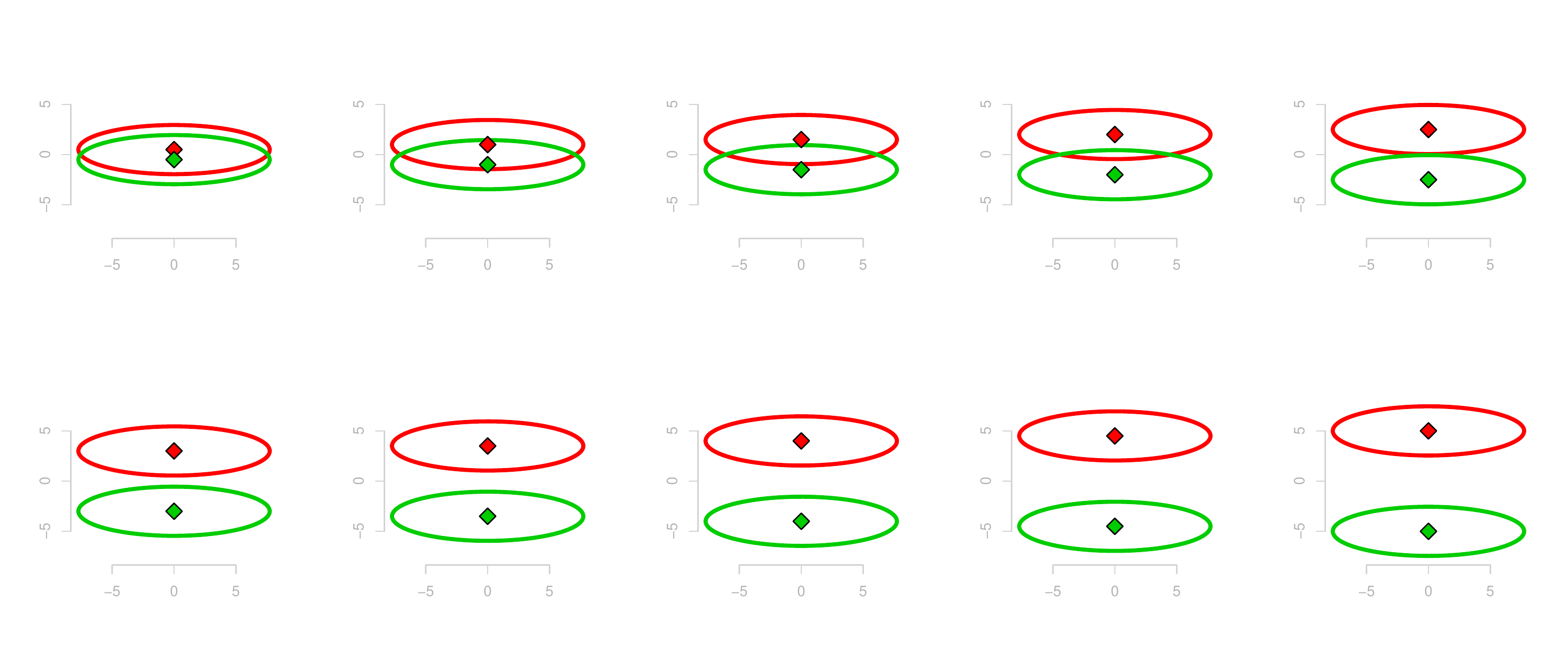}
        \caption{Diagram of increasing between cluster \textbf{distance}}\label{rys:12a}
      	 \includegraphics[width=\textwidth]{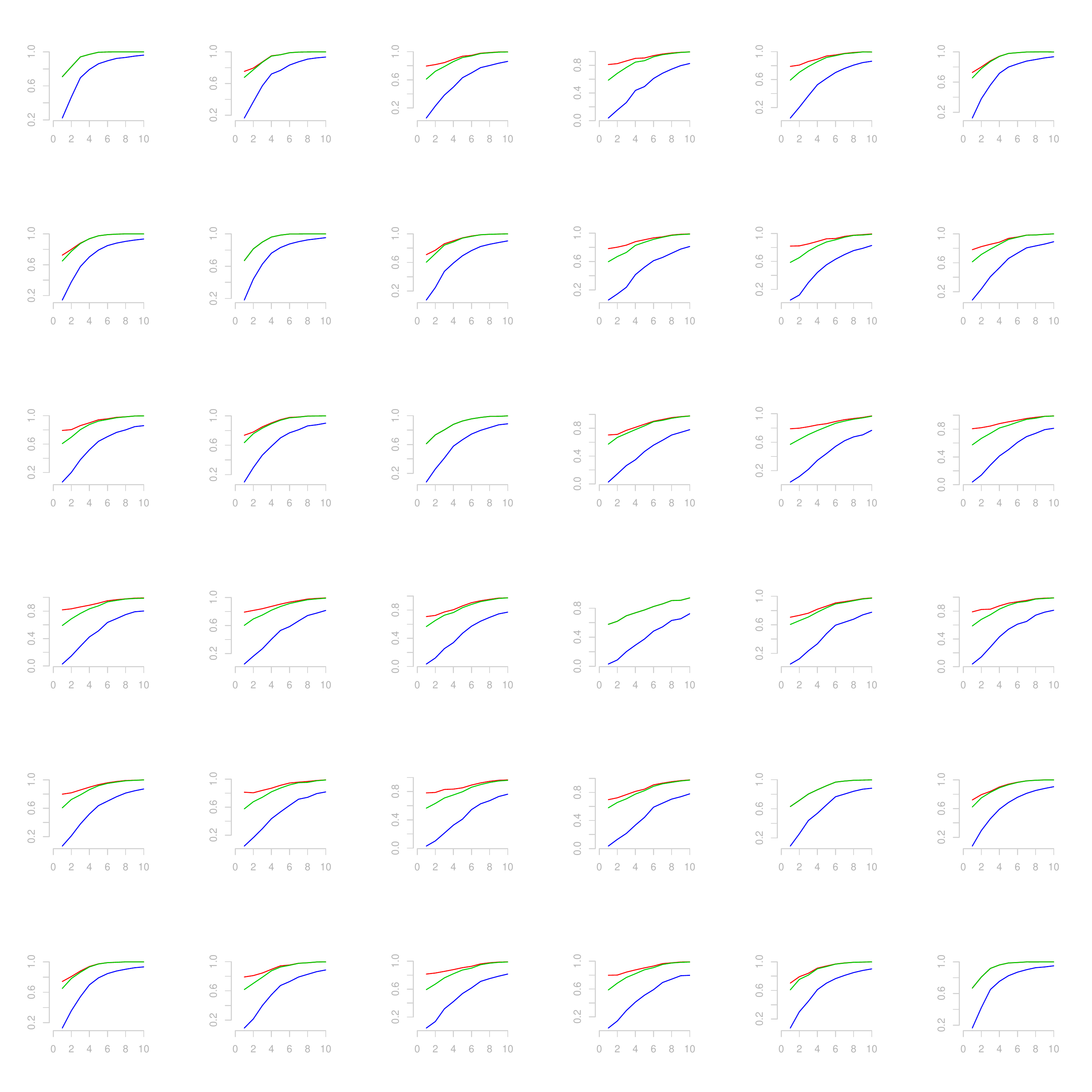}
        \caption{For clusters in position as in Figure \ref{rys:11}, each chart presents impact of increasing between cluster distance according to the pattern from Figure \ref{rys:12a}, measured with exact $(1 - \text{MLE}_{\text{err}})$ (red), its linear approximation $(1-\mathbb{P}_{\text{minmax}})$ (green) and Fisher's eigenvalue (blue).}\label{rys:12b}
   \end{center}
\end{figure}

\begin{figure}[!ht]
	\begin{center}  
       	 \includegraphics[width=\linewidth]{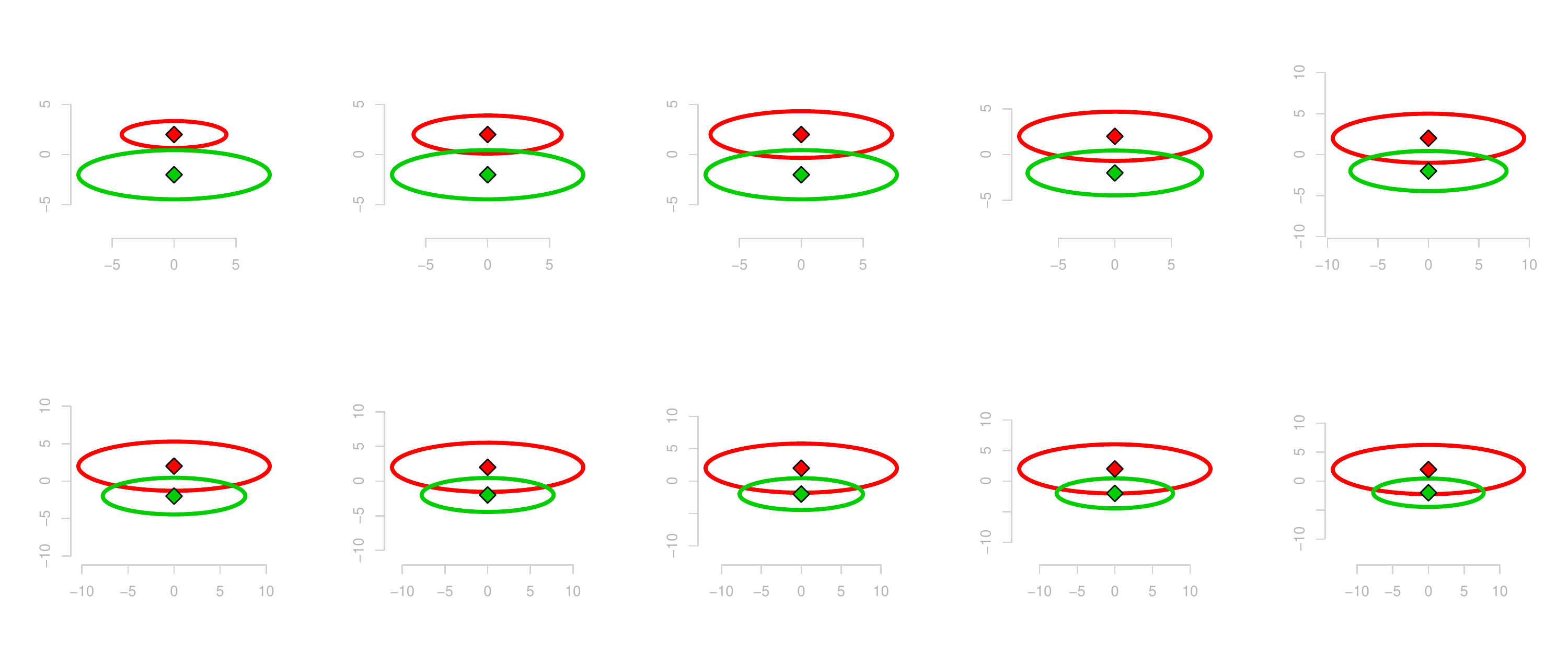}
        \caption{Diagram of increasing within cluster \textbf{dispersion -- in both spanning directions}}\label{rys:13a}
		 \includegraphics[width=\linewidth]{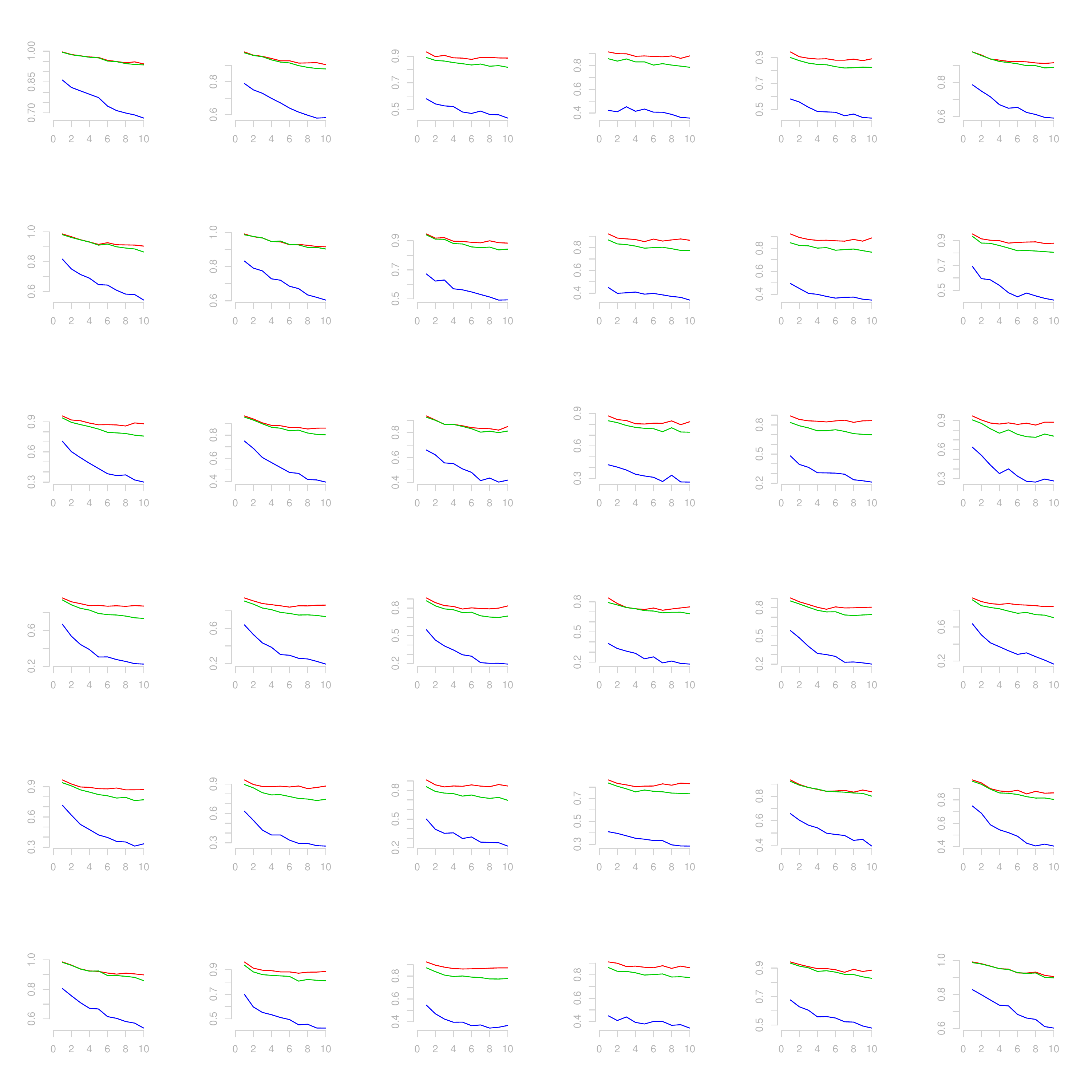}
        \caption{For clusters in position as in Figure \ref{rys:11}, each chart presents impact of increasing within cluster dispersion (both directions) according to the pattern from Figure \ref{rys:13a}, measured with exact $(1 - \text{MLE}_{\text{err}})$ (red), its linear approximation $(1-\mathbb{P}_{\text{minmax}})$ (green) and Fisher's eigenvalue (blue).}\label{rys:13b}
   \end{center}
\end{figure}

\begin{figure}[!ht]
	\begin{center}      
       	 \includegraphics[width=\linewidth]{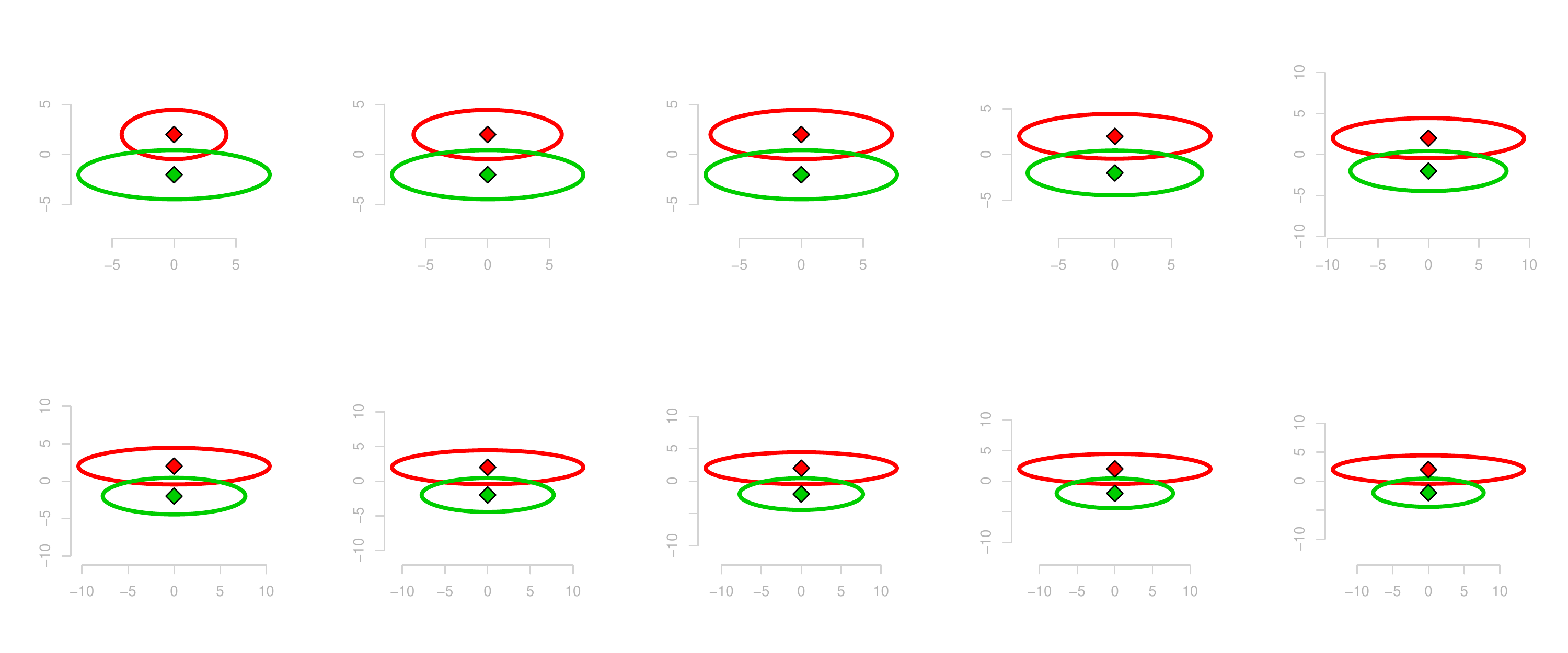}
        \caption{Diagram of increasing within cluster \textbf{dispersion -- first spanning direction}}\label{rys:14a}
         \includegraphics[width=\linewidth]{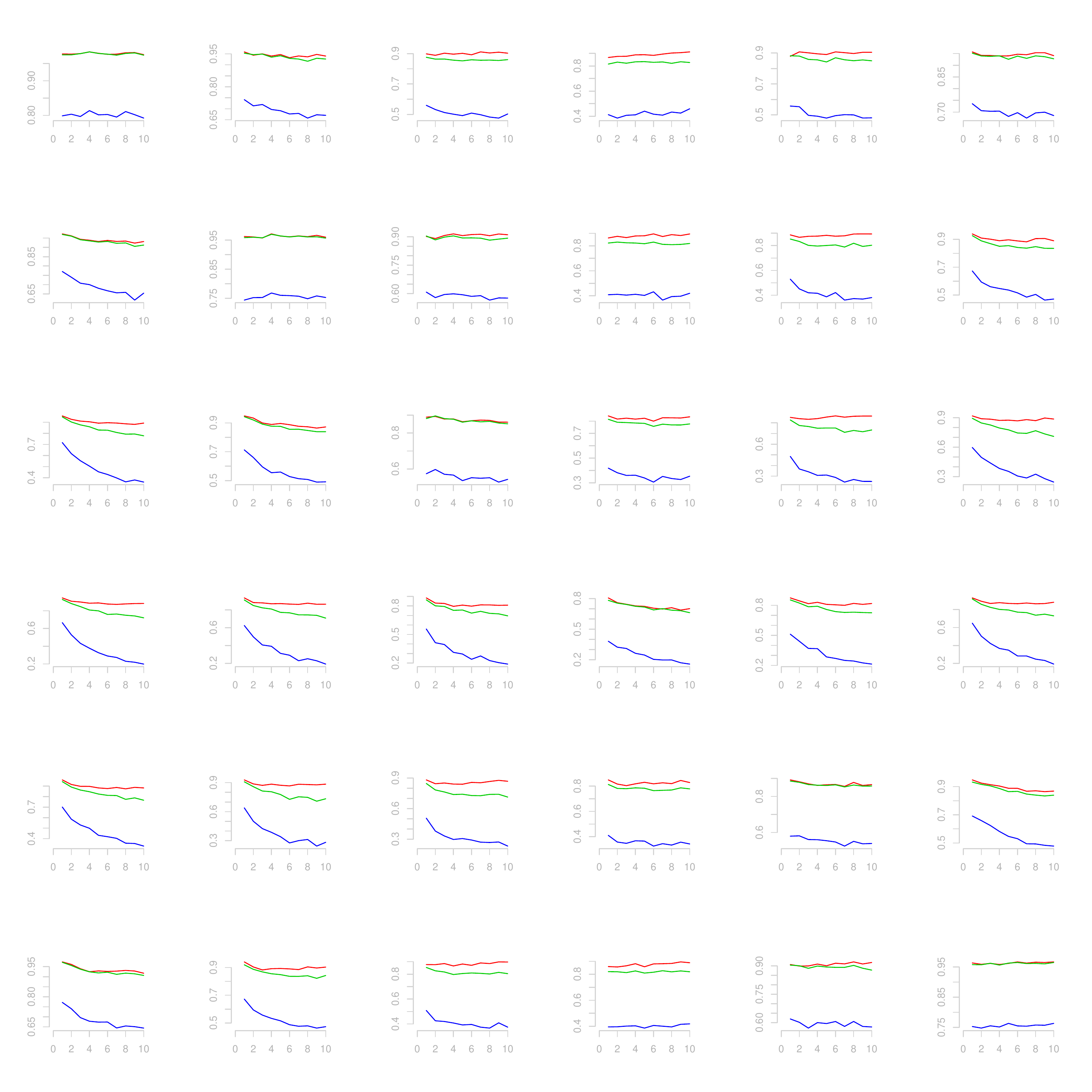}
        \caption{For clusters in position as in Figure \ref{rys:11}, each chart presents impact of increasing within cluster dispersion (first direction) according to the pattern from Figure \ref{rys:14a}, measured with exact $(1 - \text{MLE}_{\text{err}})$ (red), its linear approximation $(1-\mathbb{P}_{\text{minmax}})$ (green) and Fisher's eigenvalue (blue).}\label{rys:14b}
   \end{center}
\end{figure}

\begin{figure}[!ht]
	\begin{center}      
       	 \includegraphics[width=\linewidth, height=5.5cm]{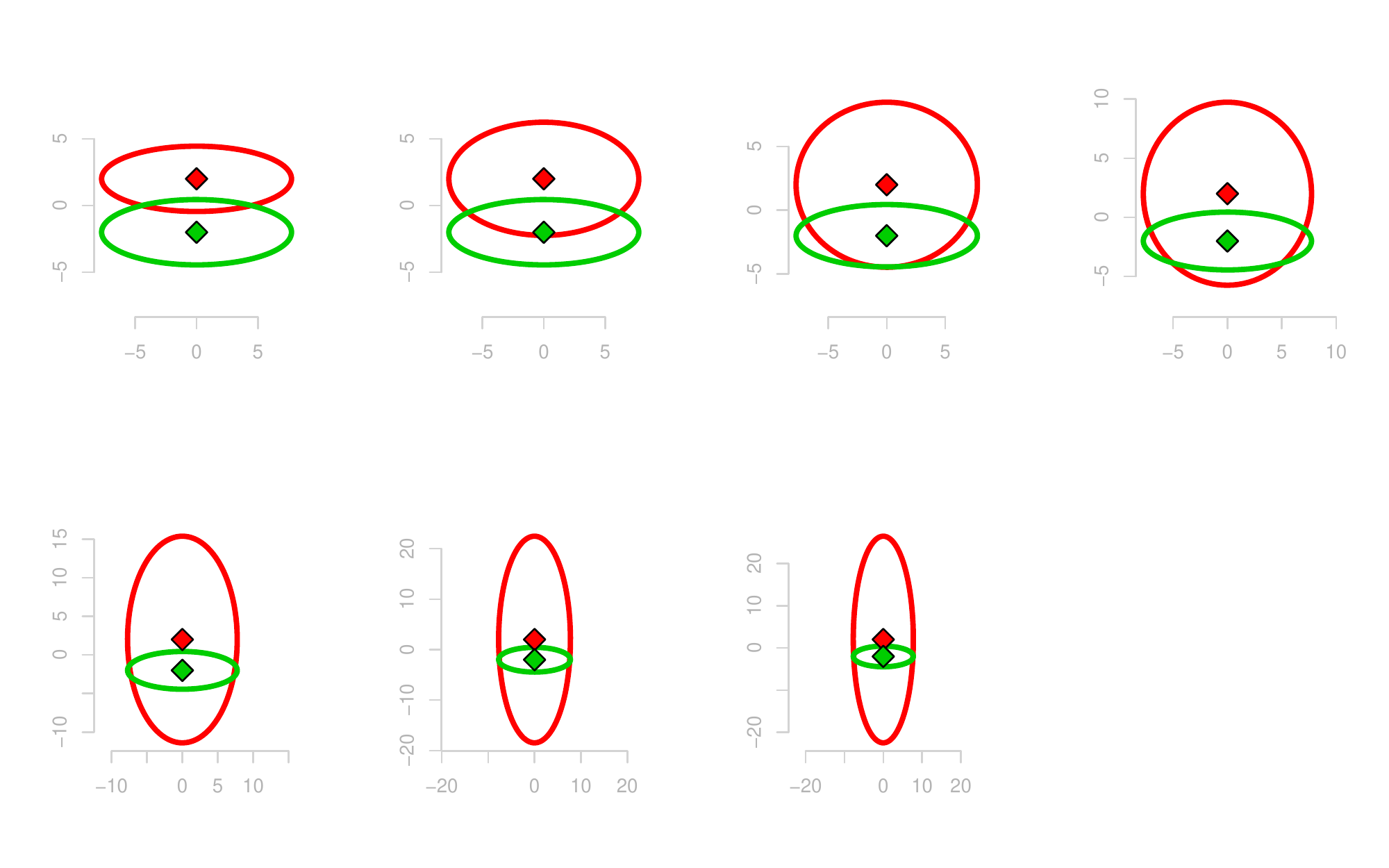}
        \caption{Diagram of increasing within cluster \textbf{dispersion -- second spanning direction}}\label{rys:15a}
         \includegraphics[width=\linewidth]{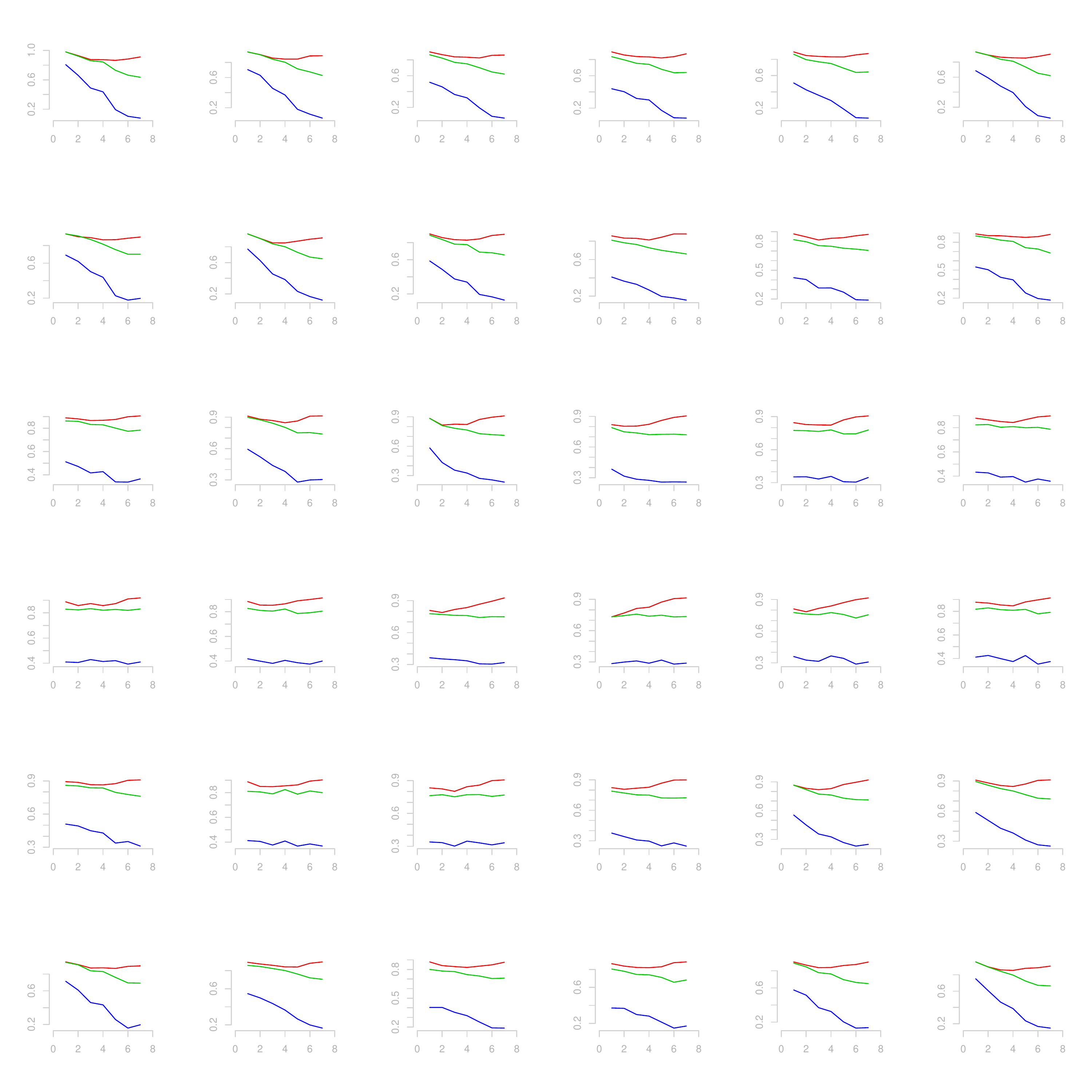}
        \caption{For clusters in position as in Figure \ref{rys:11}, each chart presents
     impact of increasing within cluster dispersion (second direction) according to the pattern from Figure \ref{rys:15a}, measured with exact $(1 - \text{MLE}_{\text{err}})$ (red), its linear approximation $(1-\mathbb{P}_{\text{minmax}})$ (green) and Fisher's eigenvalue (blue).}\label{rys:15b}
   \end{center}
\end{figure}

\begin{figure}[!ht]
	\begin{center}      
       	 \includegraphics[width=\linewidth]{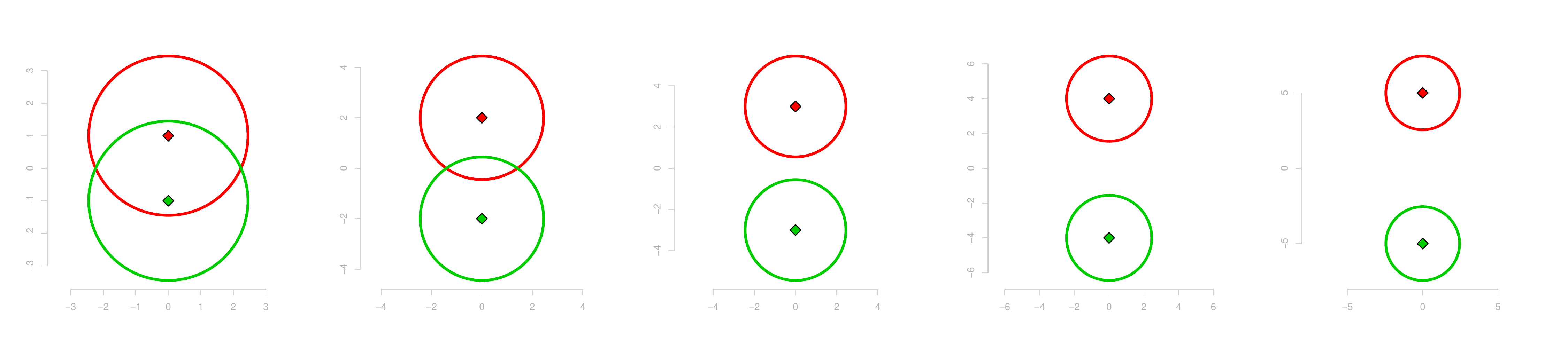}
        \caption{\textbf{Spherical components} -- diagram of increasing between cluster \textbf{distance}}\label{rys:16a}
        \includegraphics[width=\linewidth]{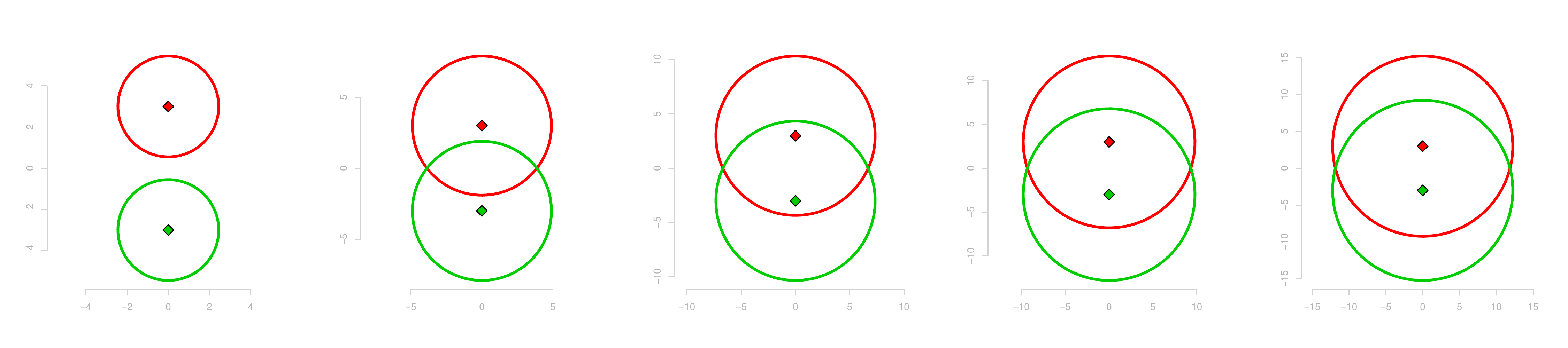}
        \caption{Diagram of \textbf{balanced} increase in within cluster \textbf{dispersion} (same increase for both clusters)}\label{rys:16b}
         \includegraphics[width=\linewidth]{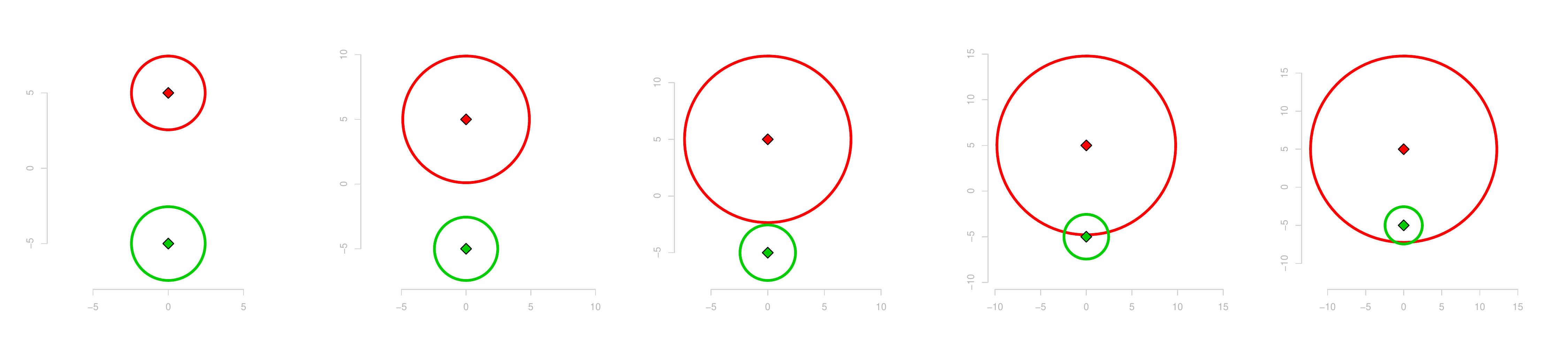}
        \caption{Diagram of \textbf{unbalanced} increase in within cluster \textbf{dispersion} (increase for one cluster only)}\label{rys:16c}
         \includegraphics[width=\linewidth, height=7cm]{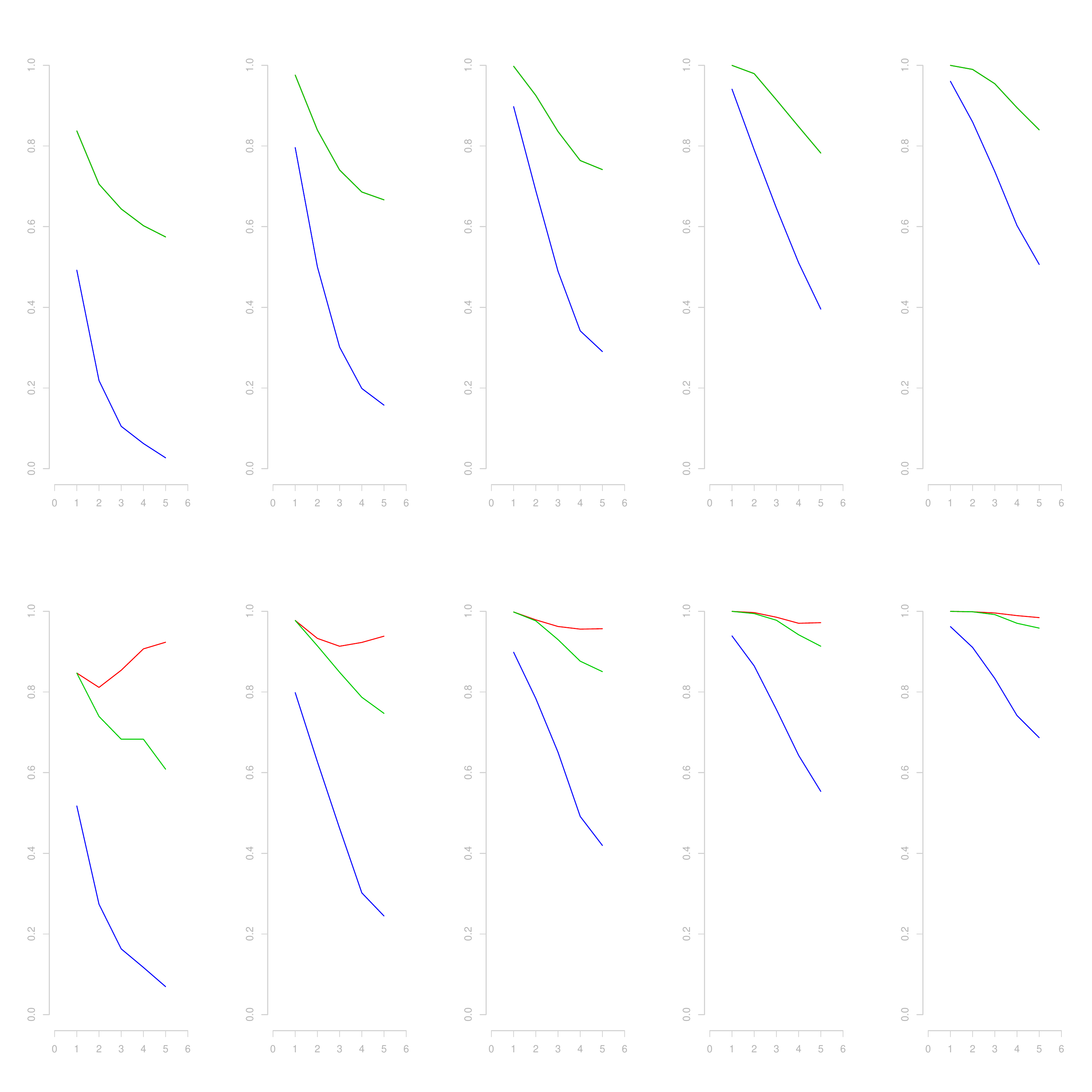}
        \caption{For spherical clusters, each line follows the distance pattern from  Figure \ref{rys:16a}, each chart presents impact of increasing within cluster dispersion -- balanced in the first line (Figure \ref{rys:16b}, unbalanced in the second (Figure \ref{rys:16c}), measured with exact $(1 - \text{MLE}_{\text{err}})$ (red), its linear approximation $(1-\mathbb{P}_{\text{minmax}})$ (green) and Fisher's eigenvalue (blue).}\label{rys:16d}
   \end{center}
\end{figure}



\begin{figure}[!ht]
	\begin{center}       
       	 \includegraphics[width=\linewidth]{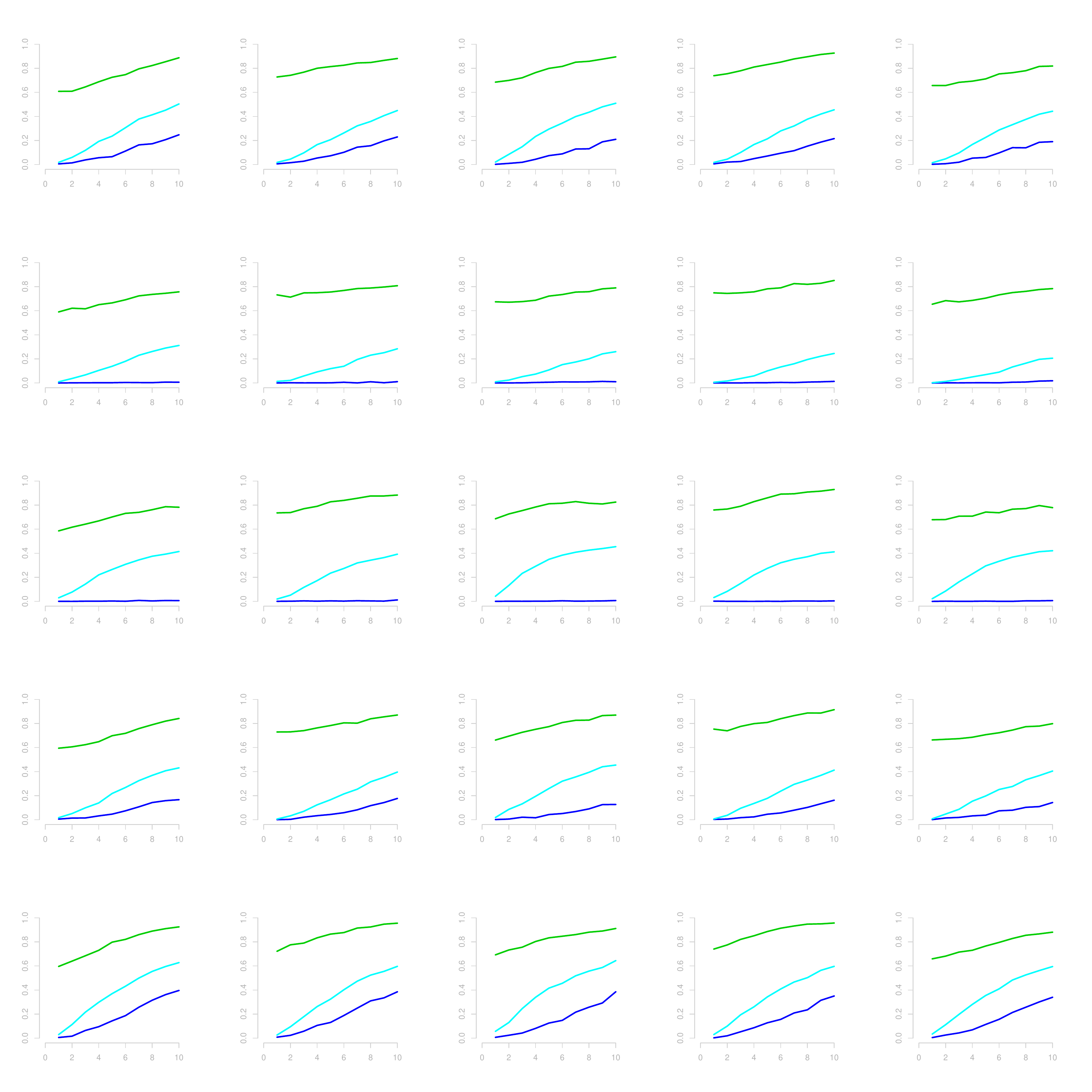}
        \caption{Three dimensions: for random (but fixed in each row) set of cluster means and random (but fixed in each column) set of covariance matrices, each chart presents impact of increasing between cluster \textbf{distance}, measured with exact integral measure (green), average Fisher's eigenvalue (turquoise) and minimum Fisher's eigenvalue (blue).}\label{rys:17a}
   \end{center}
\end{figure}

\begin{figure}[!ht]
	\begin{center}     
       	 \includegraphics[width=\linewidth]{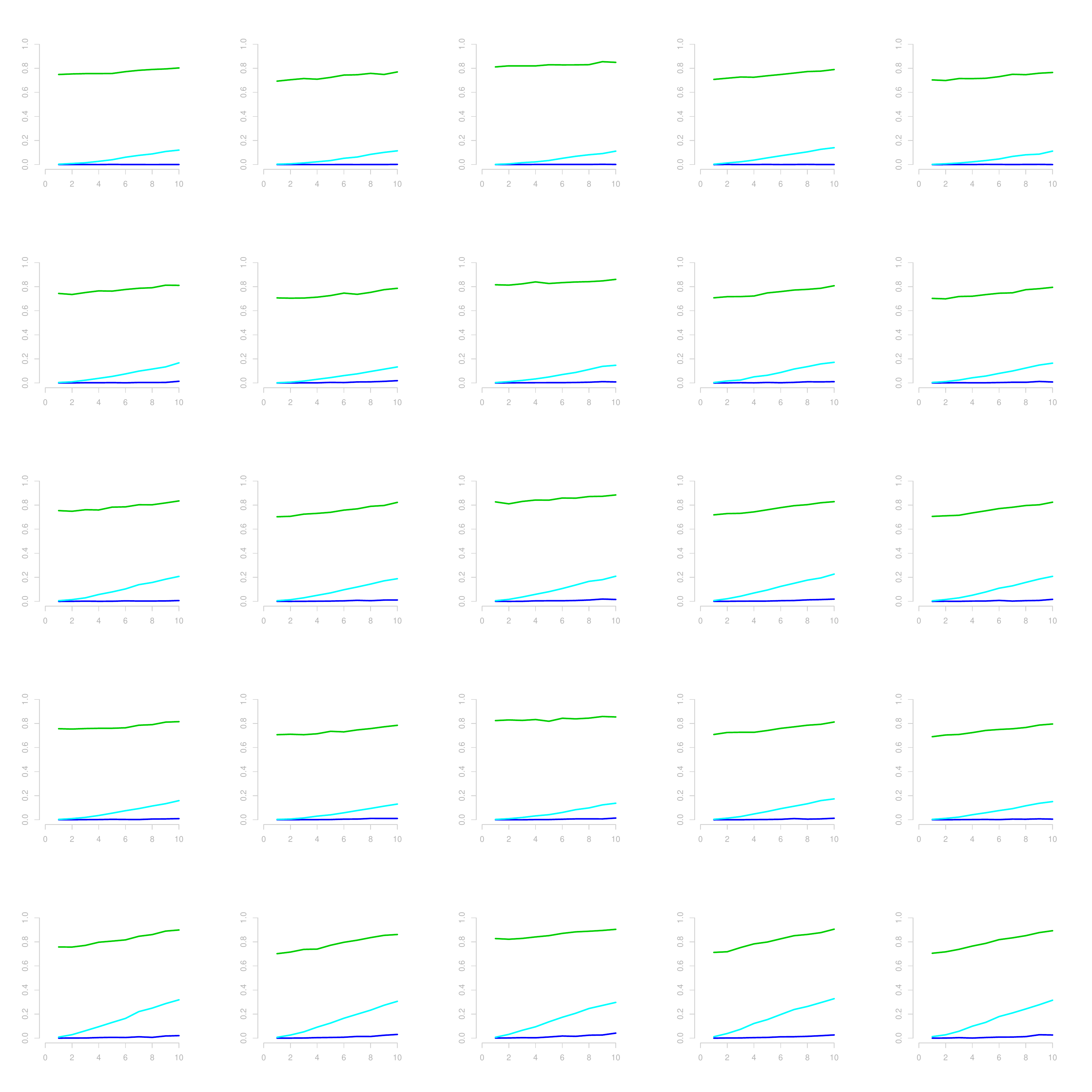}
        \caption{Five dimensions: for random (but fixed in each row) set of cluster means and random (but fixed in each column) set of covariance matrices, each chart presents impact of increasing between cluster \textbf{distance}, measured with exact integral measure (green), average Fisher's eigenvalue (turquoise) and minimum Fisher's eigenvalue (blue).}\label{rys:17b}
   \end{center}
\end{figure}


\begin{figure}[!ht]
	\begin{center}      
       	 \includegraphics[width=\linewidth]{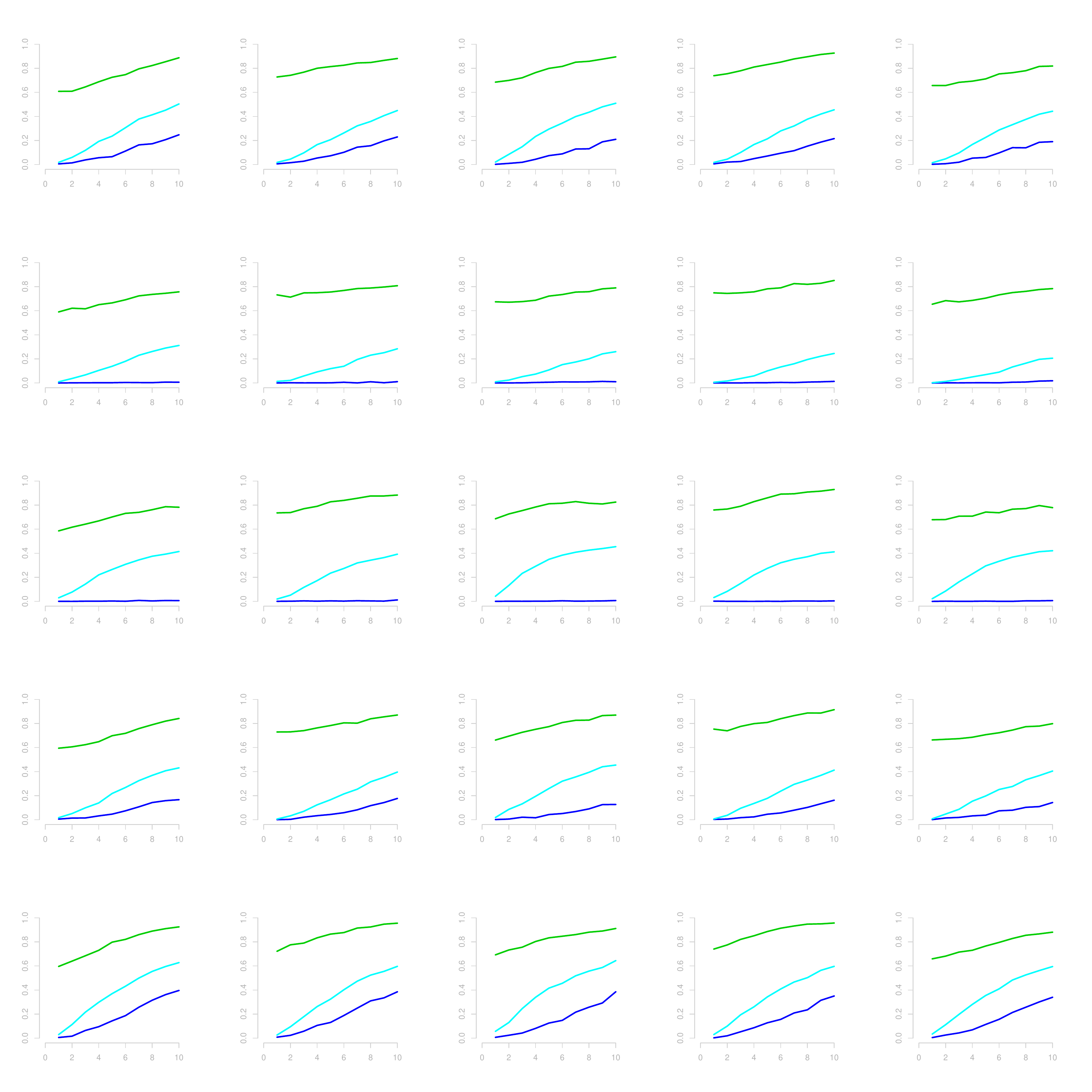}
        \caption{Three dimensions: for random (but fixed in each row) set of cluster means and random (but fixed in each column) set of covariance matrices, each chart presents impact of increasing within cluster \textbf{dispersion}, measured with exact integral measure (green), average Fisher's eigenvalue (turquoise) and minimum Fisher's eigenvalue (blue).}\label{rys:17c}
   \end{center}
\end{figure}

\begin{figure}[!ht]
	\begin{center}     
       	 \includegraphics[width=\linewidth]{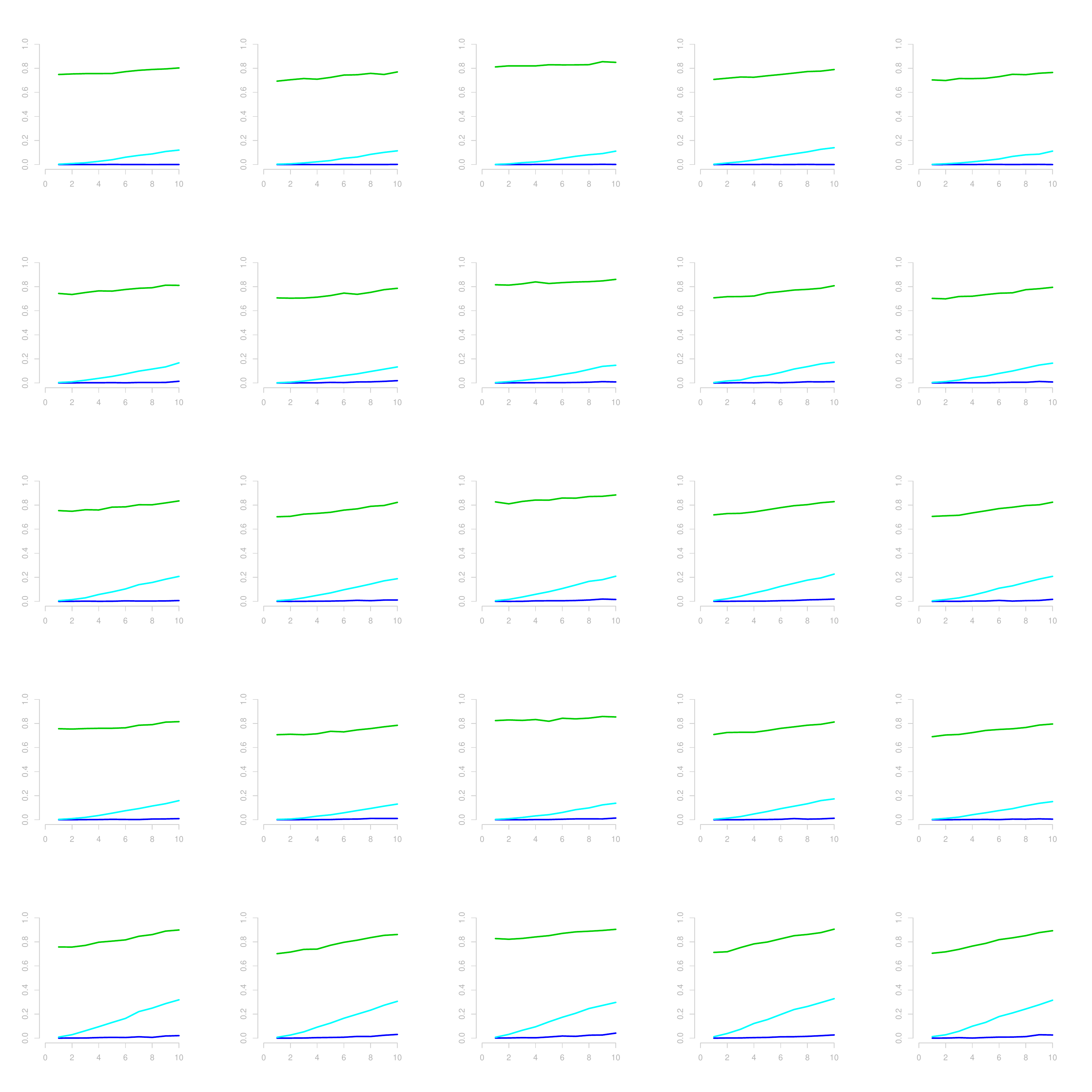}
        \caption{Five dimensions: for random (but fixed in each row) set of cluster means and random (but fixed in each column) set of covariance matrices, each chart presents impact of increasing within cluster \textbf{dispersion}, measured with exact integral measure (green), average Fisher's eigenvalue (turquoise) and minimum Fisher's eigenvalue (blue).}\label{rys:17d}
   \end{center}
\end{figure}

\end{document}